\documentclass{amsart}

\usepackage{proof}
\usepackage{geometry}
\usepackage{enumerate,xspace}
\usepackage{fancybox}
\usepackage{amsmath}
\usepackage{amssymb}
\usepackage{amsthm}
\usepackage{stmaryrd}
\usepackage[colorlinks]{hyperref}
\usepackage[all]{xy}
\usepackage{pifont}
\usepackage{tikz-cd,adjustbox,cmll}
 \usepackage{proof,latexsym,graphicx,enumerate}

\newtheorem{theorem}{Theorem}[section]
\newtheorem{lemma}[theorem]{Lemma}
\newtheorem{corollary}[theorem]{Corollary}

\theoremstyle{definition}
\newtheorem{definition}[theorem]{Definition}
\theoremstyle{definition}
\newtheorem{example}[theorem]{Example}
\theoremstyle{definition}

\theoremstyle{definition}

\theoremstyle{definition}

\theoremstyle{remark}

\theoremstyle{remark}

\newcommand{\cA}{\mathcal{A}}
\newcommand{\cB}{\mathcal{B}}
\newcommand{\cC}{\mathcal{C}}
\newcommand{\cD}{\mathcal{D}}
\newcommand{\cE}{\mathcal{E}}
\newcommand{\cF}{\mathcal{F}}

\newcommand{\cL}{\mathcal{L}}
\newcommand{\cM}{\mathcal{M}}

\newcommand{\set}{\mbox{${\bf Set}$}}
\newcommand{\rel}{\mbox{${\bf Rel}$}}
\newcommand{\Par}{\mbox{${\bf Par}$}}

\newcommand{\cat}{\mbox{${\bf Cat}$}}

\newcommand{\NN}{\mathbb{N}}

\newcommand{\stt}[1]{\stackrel{#1}{\longrightarrow}}
 
\newcommand{\Fcart}{\cF_{\!\!\!cart} } 
\newcommand{\Fmon}{\cF_{\!\!\!mon} }
\newcommand{\Fccc}{\cF_{\!\!\!ccc} }
\newcommand{\Ftop}{\cF_{\!\!\!top} }
\newcommand{\Fbtop}{{\cF_{\!\!bool}}}
\newcommand{\prim}{\ensuremath{\mathbf{Prim}}}
\newcommand{\comp}{\ensuremath{\mathbf{Comp}}}

\newcommand{\real}[2]{#1\hspace{2pt} \underline{\mathrm{\bf{r}}} \hspace{2pt} #2}
\newcommand{\eff}{\ensuremath{\mathbf{Eff}}}
\newcommand{\pass}{\ensuremath{\mathbf{Pasm}}}
\newcommand{\ass}{\ensuremath{\mathbf{Asm}}}
\newcommand{\RT}[1]{\mathbf{RT}(#1)}
\newcommand{\HEO}{\mathbf{HEO}}
\newcommand{\per}[1]{\mathbf{PER}(#1)}
\newcommand{\EC}[1]{\mathbf{EC}(#1)}

\newcommand{\scbrl}{\lbrack \! \! \lbrack}
\newcommand{\scbrr}{\rbrack \! \! \rbrack}

\newcommand{\sss}{\mathsf{s}}
\newcommand{\kkk}{\mathsf{k}}
\newcommand{\bbb}{\mathsf{b}}
\newcommand{\ccc}{\mathsf{c}}
\newcommand{\iter}{\mathsf{It}}
\newcommand{\numer}[1]{\widehat{#1}}
\newcommand{\synt}[1]{\mathsf{#1}}
\newcommand{\csynt}[1]{\mathsf{#1}}

\begin{document}

\title{Aspects of Categorical Recursion Theory}

\author{Pieter Hofstra and Philip Scott}
\address{Department of Mathematics and Statistics \\
University of Ottawa \\
Canada \\}
\email{phofstra@uottawa.ca}
\email{phil@site.uottawa.ca}
\thanks{Both authors are  partially supported by an NSERC Discovery Grant.}

\begin{abstract}
We present a survey of some developments in the general area of category-theoretic approaches
to the theory of computation, with a focus on topics and ideas particularly close to the interests of 
Jim Lambek. 
\end{abstract}

\maketitle
\tableofcontents
\setlength{\parskip}{.1cm}

\section{Introduction}

Algorithms have been discussed for thousands of years, starting with the Babylonians and later the Greeks 
(e.g. Plato's academy, Euclid in Alexandria, etc.).  These ideas were subsequently passed to (or rediscovered in) 
many mathematical cultures and civilizations (see~\cite{AngLam95}).  Indeed, the word {\em algorithm} 
itself comes from  the Latinized name of the author of a book on Hindu arithmetic, the Persian mathematician 
Muhammed ibn-M\={u}s\={a} al-Khw\={a}rizm\={i}   (c. 825).   Yet it was only in the 19th century that serious 
approaches to understanding the foundations of algorithms and computable functions began. 
For example, the modern idea of defining functions by iteration and  proofs by induction seems 
to have originated in the writings of Richard Dedekind~\cite{DedekindR:wassz}. David Hilbert's seminal lectures 
on the foundations 
of mathematics led Thoralf Skolem in the early 1920's to axiomatize the  primitive recursive functions, a  
class of inductively defined numerical functions which were intuitively computable.   Was  this all of them?  
Alas, no: a routine application of Cantor's  diagonal argument (\cite{CutlandN:comp}, p.91) shows that there are intuitively 
computable functions which are not primitive recursive.  Indeed, in 1928 Hilbert's student Wilhelm Ackermann 
constructed an explicit example of a recursively defined, intuitively computable  function which grows faster 
than any primitive recursive function.  Throughout the 1920's Hilbert  discussed the \emph{}Entscheidungsproblem 
(Decision Problem) for predicate logic, whose surprising final (negative) answer was obtained 
independently by Alonzo Church and Alan Turing in 1936, influenced by work of Kurt G\"odel (1931).  Indeed, this was the 
culmination of seminal research developing the modern theory of computability and computable functions by 
the logicians Church, G\"odel, and Church's students Stephen Kleene, J. Barkley Rosser, and 
Alan Turing in the period 1931-1936. 

Jim Lambek, in his writings and public presentations,  had a long-time interest in the foundations of 
computability and its history~\cite{AngLam95}.  His published papers include his well-known introduction of 
{\em abacuses} in 1961~\cite{Lam61a}  as a simple model of computation 
(an alternative to Turing machines) as well as his work on
 applying Gerhard Gentzen's cut-elimination algorithm and normal forms to categorical 
 coherence theory~\cite{Lam68b, Lam69, Lam72}.  
This led to his interest in typed combinatory algebras, typed lambda calculi and categorical theories of computation 
(see~\cite{Lam80} and Part III of the book with the second author~\cite{LS86}).  In linguistics and 
anthropology, 
as well as in mathematics, he often expressed interest in doing 
computation via relation algebras, e.g. from kinship terminology~\cite{BharLam95} to Mal'cev categories~\cite{CLP90} to
exact completions and partial equivalence relations~\cite{LambekJ:exampr}.  
We shall further explore some of these ideas in  Section~\ref{sec:prooftheory} below.

Here we shall examine three particular questions that occupied Lambek for many years.

  \begin{enumerate}
\item Are there {\em natural} recursion theories?
  \item What are the computable functions and functionals in various concrete categorical
structures?
  \item Are there intrinsic algebraic/categorical approaches to recursion theory?
 
 \end{enumerate}

 The  detailed discussion of these questions will be pursued in  the following sections.
As a warm up, we describe the informal meaning of Lambek's questions (and, in part, some of the associated answers).

 \subsection{On Lambek's Questions}

Concerning questions (1) and (2), in many conversations and lectures Lambek emphasized that 
``natural" recursion theories (and their classes of  computable functions)  should arise by examining 
the computable numerical functions  in various 
{\em free categories} arising in categorical logic. Here, by a \emph{free structured category} (where 
the kinds of structure we may wish to consider include monoidal structure, finite limits, cartesian closed structure, etc.)
we mean the structured category with natural numbers object (NNO) freely generated by 
the empty graph. Such a category is then initial among the categories with this structure.
We can define the notion of a representable numerical function
analogously to how this is done in mathematical logic (as in G\"odel~\cite{GodelK:ubefus}, cf. also~\cite{LS86})   and we can ask: 
which numerical functions (partial functions, functionals, etc.) are representable in the various 
relevant free categories?

To this end, these questions were taken up in Lambek and Scott (\cite{LS86}, Part III)  for categories 
associated to various higher-order logics, and will be discussed in more detail in Section~\ref{sec:compcat} below.  
We summarize some of the early literature in Figure~\ref{nat-rec}.  
  \begin{figure}[ht]
\begin{center}
\begin{tabular}{|c|c|}
\hline
{\bf Free Categories} & {\bf Definable Functions and Functionals} \\
\hline
\hline
Cartesian and monoidal with NNO & Primitive recursive functions (\cite{RomanL:carcnn, PareR:moncnn})\\
\hline
  Cartesian closed with NNO  & G\"odel's Dialectica Functionals (\cite{LS86})\\
  \hline
The free elementary topos   & Provably total functions of HAH (\cite{LS86}) and \\
with NNO &  Higher provably recursive functionals (\cite{Sced88})\\
\hline
\hline
 C-monoids and CCCs & Church's untyped lambda calculus (with surjective\\ 
 with reflexive objects & pairing) and the partial recursive functions (\cite{LS86})\\
 \hline
\end{tabular}
\caption{Natural Recursion Theories and their computable function(al)s}
\label{nat-rec}
\end{center}
\end{figure}

Concerning Question (3), since the 
1960's there has been  increasing interest in developing general categorical frameworks for 
computability theory. We mention here in particular the early work by 
Eilenberg and Elgot~\cite{EilenbergS:rec} 
on recursiveness,
and the groundbreaking work by Di Paola and Heller on recursion 
categories~\cite{DiPaolaR:domc, HellerA:exitrc, 
DiPaolaR:sompsp}, the modern incarnation of
which we shall discuss in Section~\ref{sec-turing}. Other aspects of computation and computability were studied
from a categorical standpoint by various authors; for example, Lawvere~\cite{LawvereFW:daccc} 
(see also~\cite{YanofskyN:uniasr}) gives a general 
version of the diagonal argument from which the well-known first recursion theorem, the fixed point theorem in 
untyped lambda calculus, and G\"odel's diagonalization lemma
can be obtained; while Mulry~\cite{MulryP:genbmf} introduces the \emph{recursive topos} as a natural setting to 
consider a generalization of the Banach-Mazur functionals to all higher types.

We shall not attempt to give a full historical account of the large recent literature on categorical recursion theory. 
Instead, we shall focus on a few areas that the authors have become involved in which represent new 
directions of independent interest, but which also overlap with Lambek's interests.

\subsection{General Notation and Background}
We now introduce some notation and terminology for some of the structure appearing frequently in this paper.
We assume that the reader is familiar with basic category theory.
Standard references include~\cite{MacLaneS:catwm, AwodeyS:catt}. 
Some familiarity with lambda calculus~\cite{BarendregtH:lamcss} and
the basic theory of computation~\cite{CutlandN:comp, OdifreddiP:clart} is also an advantage. 

The (large) category of sets and functions is denoted by $\set$. The category of sets and partial functions 
(that is, single-valued relations) is denoted by $\Par$. By $\rel$, we mean the category of sets and relations. 

An \emph{idempotent} is an endomorphism
$e:A \to A$\/ for which $ee=e$. An idempotent $e$\/ \emph{splits} when there is an embedding-retraction pair
$m:B \to A, r:A \to B$\/ with $rm=1_B, mr=e$. 
The \emph{idempotent splitting} or \emph{Karoubi envelope} of $\cC$\/ is
the category $\mathcal{K}(\cC)$\/ whose objects are the idempotents of $e$, and whose morphisms $(A,e) \to (B,d)$\/ are
maps $f:A \to B$\/ with $df=f=fe$. More generally, for $E$\/ a set of idempotents $,\mathcal{K}_E(\cC)$\/ is the
full subcategory of $\mathcal{K}(\cC)$\/ on the objects determined by $E$. 
Finally, a \emph{retract} of an object $A$\/ is an object $B$\/ together with an 
embedding-retraction pair $m:B \to A, r:A \to B$, with $rm=1_B$. We write $B \triangleleft A$\/ to indicate that $B$\/ is a 
retract of $A$.

\section{What is a computable function?}\label{sec:computable}
As explained in the introduction, 
one of the fundamental challenges for mathematical logic in the first quarter of the 20th century was to make
precise the notions of \emph{computation}, \emph{computable function}, and \emph{computable set}.\footnote{Originally,
the terminology \emph{recursive function} has been prevalent, due to the emphasis on the use of recursive procedures.
Indeed, many authors have referred to the subject as \emph{recursion theory}.
However, as argued in~\cite{SoareR:comr}, the term \emph{computability} more aptly captures the flavour of the subject,
and also emphasizes inclusion of other notions or models of computation.} 
This section reviews some of these developments of the classical theory, setting the stage for the categorical
approaches to be introduced later.

\subsection{The primitive recursive functions}\label{subsec:primrec}
As summarized in the Introduction, historically the attempts to define computability focussed on iterative or 
recursive procedures.  These seem to have been first analyzed in the writings
of Dedekind in the 19th century~\cite{DedekindR:wassz}.  
A formal system of  Primitive Recursive Arithmetic, concomitant with 
Hilbert's foundational lectures in the 1920's, 
was developed by Skolem~\cite{SkolemT:fouea}. R\'ozsa P\'eter's 
work in the early 1930s (later presented as~\cite{PeterR:recfns})
is considered to have provided the foundations for the theory of recursive functions;
many of the central ideas were 
further developed in detail by Hilbert-Bernays~\cite{HilbertD:grum}, and especially 
Goodstein~\cite{GoodsteinR:RNT}.  Moreover,  these functions 
 were also used by G\"odel in his famous Incompleteness Theorem paper~\cite{GodelK:ubefus}.

Consider total numerical functions $\NN^k\rightarrow \NN$, $k\geq 1$.  We recall
the traditional definition of primitive recursion, then include a somewhat non-standard definition by Lambek.

 \begin{definition}[Primitive Recursive Functions]\rm 
 The {\em primitive recursive functions} are the smallest class $ \prim$ of numerical functions 
 generated from Basic Functions by composition (or substitution) and primitive recursion.  

The Basic Functions are the constant zero function $Z(x)=0$, the successor function $S(x)=x+1$, and the 
projection functions $U^n_i:  \NN^n\rightarrow \NN, \; U^n_i(x_1, \ldots, x_n)=x_i$.  The closure rules are as follows
($\vec{x}$\/ denotes an element of $\NN^n$):

\noindent
\begin{itemize}
\item {\bf Composition}:   if $f_1: \NN^n\rightarrow \NN \in \prim$ and
 $g(u_1, \cdots , u_k): \NN^k\rightarrow  \NN \in \prim$, then 
 $comp(g, f_1,\ldots, f_k): \NN^n\rightarrow \NN  \in \prim$,
 where 
\[
 comp(g, f_1,...,f_k)(\vec{x}) = g(f_1(\vec{x}),\ldots, f_k(\vec{x}))   .
\]
\item {\bf Primitive Recursion}:
 if $g(\vec{x}), h(\vec{x}, y, u)\in \prim$ then so is $rec(\vec{x}, y)$, where  
 \[ rec(\vec{x}, 0) = g(\vec{x}), \ \  rec(\vec{x}, S(y)) = h(\vec{x}, y, rec(\vec{x}, y)).\]
 \end{itemize}
 \end{definition}
 
Using this, we may also define a relation $R \subseteq \mathbb{N}^k$\/ to be primitive recursive when
its characteristic function is. Most of the numerical functions and relations used in everyday mathematics
are primitive recursive. 


 Let us mention a somewhat nonstandard definition\footnote{We modify slightly
 the Basic functions, which were missing one.} of primitive recursive functions, 
 introduced by Lambek in~\cite{AngLam95}, p. 246. Following Lambek, we elide function 
 arguments, writing e.g. 
  $f\vec{x}u\vec{v}$ for $f(\vec{x},u,\vec{v})$, etc.
  \begin{definition}[Lambek's Primitive Recursive Functions]
  \begin{enumerate}[(i)]
  \item Basic functions:  {\em Identity} $Ix = x$,  {\em Successor}  $Sx=x+1$, and {\em Zero} $Zx = 0$.
  \item Generating Rules:
  \begin{enumerate}[(a)]
  \item {\em Substitution}:   given functions $f\vec{x}u\vec{z}$ and $g\vec{y}$
  we can form  $h\vec{x}\vec{y}\vec{z} = f\vec{x}(g{\vec y})\vec{z} $.
  \item {\em Interchanging two arguments}:  given $f\vec{x}uv\vec{y}$, we can form
  $g\vec{x}uv\vec{y} =  f\vec{x}vu\vec{y}$.  
  \item {\em Contracting two arguments}: given $f\vec{x}uv\vec{y}$, we can form
  $g\vec{x}u\vec{y}= f\vec{x}uu\vec{y}$.
  \item {\em Introducing dummy arguments}: given  $f\vec{x}\vec{y}$, we can form
  $g\vec{x}u\vec{y}= f\vec{x}\vec{y}$.
  \item {\em Primitive Recursion}:  given $g\vec{x}$ and $h\vec{x}yz$, we can form
  $f\vec{x}y$, where $f\vec{x}0 = g\vec{x}$,  $f\vec{x}(Sy) = h\vec{x}y(f\vec{x}y)$.  \end{enumerate}
 \end{enumerate}
  \end{definition}
 Note that Lambek's rules generating $\prim$ are closely related to the Curry-Howard 
 functional interpretation of intuitionistic sequent calculus proofs (with non-logical axioms).  Indeed, consider a proof
 of an intuitionistic sequent $A_1,\cdots,A_n \vdash B$.  The  {\em functional interpretation}  
 interprets the proof by functional ``proof terms" (see Girard~\cite{GLT89} and Lambek~\cite{Lam87})
of the form
\[ x_1: A_1,\cdots, x_n: A_n \vdash f\vec{x} :B.\] 
 The identity function interprets
 the identity sequent $A\vdash A$.  The generating rules (a)-(d) above correspond respectively to 
 interpreting the following rules of sequent calculus (by associating to proof terms for each of the premises
 a proof term of the conclusion): cut (a), interchange (b),
 contraction (c) and weakening (d). 
 \[
\begin{array}{cccc}
\infer[\textrm{cut}]{\Gamma,  \Sigma , \Delta \vdash C}{\Gamma, A,  \Delta \vdash C & \Sigma \vdash A} & ~~ & ~~ &
 \infer[\textrm{interchange}]{\Gamma, B, A , \Delta \vdash C}{\Gamma, A, B, \Delta \vdash C} \\
 \\
 \infer[\textrm{contraction}]{\Gamma, A,\Delta \vdash B}{\Gamma, A, A, \Delta \vdash B} &&&
  \infer[\textrm{weakening}]{\Gamma, A,\Delta \vdash B}{\Gamma, \Delta \vdash B} \\
\end{array}
\] 
 Finally, the zero function, successor, and primitive recursion (e) above may be thought-of as
 non-logical axioms or rules specifying a weak natural numbers object, i.e. a particular
 type $\synt{N}$.   In this sense, we have distinguished proofs
 $\vdash 0: \synt{N}$ and $x:\synt{N}\vdash Sx: \synt{N}$, and primitive recursion is a special case
 of the iterator (see Subsection \ref{subsec:NNOs} below).
\subsection{The computable functions}
So, does Computable = Primitive Recursive?  Alas, no, by a standard application of  Cantor's diagonal argument.  
Indeed, the previously mentioned Ackermann function (which is computable but not primitive recursive) 
can be defined by a so-called \emph{double recursion scheme} (see e.g. R. P\'eter's book~\cite{PeterR:recfns}).  
So what is a computable function?  
This was taken up in a remarkable development in the years 1931--1937  (primarily centered around Princeton University)
which, as it turned out,  led to the foundations of modern computer science.
Let us briefly recall the history.

\begin{itemize}
\item  A. Church (1932-34) and his students (S. C. Kleene, J. B. Rosser)  developed (untyped) lambda calculus as 
a model of computation (and, as later realized in the 1960's,  a foundation of modern programming language theory).    
Church formulated {\em Church's Thesis (1936)}: the intuitively computable numerical 
functions are exactly those you can compute in $\lambda$-calculus. This thus answered the age-old question we
began with.

  \quad Originally, however, Church's thesis was not believed by G\"odel (there being insufficient evidence at the time).  
   However,  in rapid developments, new evidence arose:

\item Kleene (1934-35) developed the partial $\mu$-recursive functions:  
we add to $\prim$ the following generating scheme on partial functions, called   {\em minimalisation}: given $g(\vec{x},y)$, 
we can form 
\[
f(\vec{x}) = \mu y. g(\vec{x},y) = 0
\] 
where $\mu y. g(\vec{x},y) = 0$ means {\em the  least $y$ such that $g(\vec{x},y) = 0$} \footnote{Provided for all $z < y,
g(\vec{x},z)$ is defined and not = 0. If there is no such $y$, $\mu y. g(\vec{x},y) = 0$ is undefined. } 
(If we wish to restrict to total functions, we add the proviso $\forall \vec{x} \exists y . g(\vec{x},y) = 0$.)
\item G\"odel-Herbrand (1934).  G\"odel lectured on an equation calculus to define ``computable" functions, 
based in part on a letter from Herbrand.  This is described
in Kleene's book~\cite{KleeneS:intmm}.
\item Turing (1936) independently introduced Turing machines:  an abstract mechanical computing device. 
He gave a convincing analysis of the meaning of  being ``computable"  without restrictions on space or time.  
This led to {\em Turing's thesis}: the intuitively computable functions were those computable by Turing's abstract machines. 
This ground-breaking paper also showed the recursive unsolvability of Hilbert's Entscheidungsproblem,
simultaneously and independently solved by Church in 1936  (who was inspired by his studies in untyped lambda calculus).


\item  Turing then became a student of Church at Princeton. During the period
1936-37,  Church, Kleene and Turing  carefully proved the ``equivalence" of the above different models of computability, 
in the sense that all notions gave exactly the same class of computable functions!  
This work convinced G\"odel of the truth of the Church-Turing thesis (CT).\footnote{CT is not a mathematical statement:   
it is an experimental statement,  identifying an informal class 
(namely, the ``intuitively computable" numerical functions) with a precise mathematical class of functions. }
\end{itemize}

\subsection{Some Newer Models of Computability}

After the exciting results in the late 1930's, mathematicians continued the analysis of abstract theories of computing.  
For example   Emil Post (1943) and Andrei Markov (1951)
developed theories of  computability based on string rewriting grammars  (following in the
footsteps of the Norwegian mathematician Axel Thue).  These notions of computability turned out to be Turing complete, i.e.,  
equivalent to Turing computability. In 1944, Post~\cite{PostE:recesd} also initiated the systematic study 
of the recursively enumerable sets (previously defined
by Kleene and Church in terms of images of recursive functions), in particular the study of the r.e. degrees. 

A particularly interesting period in the more recent modelling arose in 1960-61   (simultaneously and almost independently): 
the development of  {\em Unlimited Register Machines}.  Within a period of a few months, papers by  
J. Lambek,  Z. Melzak, M. Minsky, and (slightly delayed) J. Shepherdson and J. Sturgis introduced 
this influential model of computability.\footnote{Lambek's and Melzak's papers appeared back to back
in the same issue of the Bulletin of the Canadian Mathematical Society. Lambek's paper
is a considerable simplification of Melzak's approach.}     

Lambek's paper~\cite{Lam61a} was by far the  simplest to read of all the papers on Register Machines, 
and used a highly graphical syntax, akin to flowcharts.  Register machines were particularly influential pedagogically, 
compared to the intricacies of Turing machines.  A direct translation between Lambek's machine models
 and Turing Machines is given
in Boolos and Jeffrey~\cite{BoolosG:coml}.

 Let us briefly recall the formalism.  A \emph{Lambek abacus} consists of a series of {\bf Locations}
 (or registers) of arbitrary capacity (denoted   $X$, $Y$, $Z$, $\cdots$ ), into which we
  may put (or remove) pebbles, called {\bf Counters}.  We assume an unlimited supply of (indistinguishable) 
  pebbles as counters.  There are a small number of 
 {\bf Elementary Instructions} for building abacuses, as follows:
 
 \begin{figure}[ht]
 \begin{center}
 $
 \begin{array}{c}
{\mbox{\sf \bf Start}}\\
 \downarrow\\
 \end{array}
 $
 \qquad\qquad$
 \begin{array}{c}
 \downarrow\\
{\mbox{\sf \bf Stop}}\\
 \end{array}
 $\qquad \qquad
$
 \begin{array}{l}
 \downarrow\\
\!X^+\\
 \downarrow\\
 \end{array}
 $
 \qquad  \qquad
 $
 \begin{array}{rcl}
 &\downarrow& \\
&Y^-&\\
 \swarrow\!\!\!\!&&\!\!\!\!\searrow\\
 \end{array}
 $
 \ \ \begin{minipage}[c]{1in}
 (If $Y\not = \emptyset$, take one pebble away and go to the left; else go to the right)
\end{minipage}
\caption{Abacus instructions}
\label{fig:abacus}
\end{center}
\end{figure}

\noindent
Here $X^+$ denotes the operation of adding one pebble to location $X$.
Programs are formed from a finite number of instructions, arranged in a flow chart (directed graph) with root {\bf Start}, 
possibly with {\em feedback loops.} In section~\ref{subsec:trace} we will discuss the categorical semantics of such
a graphical notion of computation.

 \section{Lambek's Categorical Proof Theory}\label{sec:prooftheory}
 
 Categorical logic is concerned with the study of classes of categories with additional categorical 
 structure, such as categories with finite limits,
 regular categories,  monoidal (closed) categories, cartesian closed categories, first-order categories, toposes, and so on.
 Ideally, such a class of categories corresponds to a well-behaved fragment of logic; for example, cartesian closed categories
 correspond to typed lambda calculus (see below). This correspondence means that there is a sound and complete 
 interpretation of the logic
 in this class of categories. On the one hand, this allows us to use proof-theoretic techniques
 (rewriting for example) to reason about categorical structure, while on the other hand we may apply categorical
 results to obtain information about logical systems. Categorical proof theory is particularly concerned with
 the study of syntactically generated categories and their properties. This section describes some of the
contributions due to Lambek, as well as some related developments. 
 
 \subsection{A brief history}\label{subsec:proofhistory}
 Lambek's  early works in mathematical linguistics~\cite{Lam58, Lam61b} 
 as well as  his later work in categorical coherence theory~\cite{Lam68b,Lam69,Lam72}
employed proof theory, notably Gentzen's sequent calculi. 
Coherence theorems in category theory were aimed at answering the following very general question: 
(see Mac Lane~\cite{MacLaneS:catwm})  {\em given a freely generated structured category $\cC$, prove that 
every diagram (built from some canonical morphisms) commutes.}   
Lambek  reformulated the question more generally as follows:
\begin{enumerate}[(i)]
\item Given a freely generated structured category $\cC$, how do we effectively generate the hom-sets   $Hom_{\cC}(A, B)$?
\item Find an effective method to solve the word problem for hom-sets in such $\cC$.  In particular,  any two
morphisms with the same domain and codomain generated from the canonical morphisms must be equal.
 \end{enumerate}

Lambek's seminal idea was to reformulate this problem using proof-theory, then apply Gentzen's Cut-Elimination 
(or Normalization) theorems. Namely, he considered freely generated monoidal or residuated categories as kinds 
of ``logics" or ``labelled deductive systems":   the objects of such categories are ``formulas" 
(freely generated from some atomic ones), while arrows would then be {\em equivalence classes of} proofs (or proof trees). 
 
 In particular, an arrow $f\in Hom_{\cC}(A,B)$ would be considered as a proof of the Gentzen sequent  $A\vdash B$, 
 while composition of arrows  $f:A\to B$ and \\ $g:B\to C$ to obtain $g \circ f: A\to C$ becomes an  instance of the Cut-Rule.  
 The equations of a category force one to  impose  the notion of ``equality of proofs". 
 Algebraically, one generates a {\em congruence relation} on proofs (or better, between proof trees).

 For (i), we generate all proofs of the sequents $A\vdash B$, by Gentzen's {\em proof search}.  
 For the word problem (ii), Gentzen's cut-elimination methods amount to introducing a compatible 
 rewriting system on proofs.  To decide if two proof trees denote the same arrow or not,   
 reduce each to a unique {\em normal} (or cut-free) form.  The problem of deciding equality  of arrows 
 amounts  to deciding if their normal forms are identical or not.\footnote{An equivalent formulation~\cite{TroelstraA:baspt} 
 of a coherence theorem for a free category of some kind says:  
 given any two objects $A,B$, there is at most one  proof (built from canonical arrows) of the associated sequent $A\vdash B$.}

 Lambek pursued these ideas in the late 1960's and early 1970's using cut-elimination to solve the word problem 
  for (among others) residuated and biclosed monoidal   categories in~\cite{Lam68b,Lam69,Lam72}.  
  But  it was soon realized by proof theorists, beginning with G. Mints~\cite{Mints:CCTP},  
  that natural deduction calculi (and their associated lambda calculi of proof-terms, under normalization) 
  leads to a smoother technical framework for such word problems. 
  Mints and his students greatly increased the scope of Lambek's proof-theoretic approaches to coherence, 
  influencing even Kelly and Mac Lane~\cite{KellyMac71}.  
  Normalization approaches to coherence/decision problems for monoidal categories 
  (using reduction of lambda-like proof terms) were first investigated by Mints and his students (\cite{Mints:CCTP,Mints:PTCT},
  reprinted in~\cite{Mints:select}).  In the case of  monoidal closed categories, it was shown in Mac Lane~\cite{MacL82} that 
  Mints' proof-theoretic methods agreed almost exactly with the approach to coherence due to Kelly and Mac Lane, 
  all of which in turn were influenced by Lambek's original use of Cut-Elimination.

Meanwhile, in the 70s and 80s, 
Lambek's own algebraic studies on functional completeness and combinatory logics~\cite{Lam74,Lam80}, 
led him to consider connections of lambda calculi to
freely generated cartesian and cartesian closed categories.
Around the same time, work in computer science in applying lambda calculi and 
natural deduction to functional languages   led to the now-common practice of 
assigning lambda- (or proof-) terms to proof trees~\cite{GLT89}. Hence ``equality of proofs" becomes provable equality of the associated  terms assigned to the proof trees.  This is sometimes
known as the {\em Curry-Howard-Lambek correspondence}, to be discussed in more detail below.


After the introduction of Girard's Linear Logic in 
1986~\cite{Gir87} (which used sequent calculi and gives particular analysis of the 
structural rules)  Lambek realized his earlier work in linguistics amounted to a kind
of substructural (linear) logic without structural rules.  He introduced generalizations of 
deductive systems to more general Gentzen sequents with their associated {\em multicategories} and  
term calculi~\cite{Lam87}. On the subject of categorical proof theory, cut-elimination and applications to 
(structured) monoidal categories, linear logics, coherence theorems, et cetera, 
there has been an explosion of activity.  As a small sample of the extensive literature, we mention works of  R. Blute, R. Cockett, R. Seely
and co-workers~\cite{BCST96,CocSeel97, BCS97}
K. Do\u{s}en, et al.~\cite{Dos99,DosPet04,DosPet07}, B. Jay~\cite{Jay89,Jay90}.

\subsection{Internal Languages and free categories}
As mentioned above, coherence problems are often formulated in terms of free categories.  Let us make this
more precise. Suppose that $\mbox{\bf S-Cat}$\/ is a category whose objects are structured categories
and whose morphisms are structure-preserving functors. There is a forgetful functor
\[ \mbox{\bf S-Cat} \longrightarrow \mbox{\bf DirGrph} \]
to the category of directed graphs.
The free structured category generated by a (small) graph $\mathbb{G}$, denoted $\cF(\mathbb{G})$,  can be 
described in terms of a left adjoint to this forgetful functor. In~\cite{LS86} this left adjoint is constructed using logical syntax
along the following lines. 
\begin{enumerate}[(i)]
\item  One sets up an equivalence of categories
$
\mbox{\bf S-Cat}\stt{\simeq} \mbox{\bf Lang}
$
where \mbox{\bf Lang} is 
some category of formal theories (whose morphisms are ``interpretations"  which 
preserve the structure exactly).  The equivalence is implemented by a pair of functors:
$L: \mbox{\bf S-Cat}\longrightarrow \mbox{\bf Lang}$ which associates to every category
$\cC$ a so-called {\em internal language} and $C: \mbox{\bf Lang}\longrightarrow \mbox{\bf S-Cat}$, 
which associates to a language $\cL$, a category
$C(\cL)$, called the (syntactic) category {\em generated by $\cL$}.  
\item Next, one constructs, given a directed graph $\mathbb{G}$, the theory $\cL_{\mathbb{G}}$
generated by $\mathbb{G}$.
The types of $\cL_{\mathbb{G}}$\/ are generated from the nodes of $\mathbb{G}$, while the terms
are generated using the term-formation rules of the logic by including 
the arrows of ${\mathbb{G}}$\/ as term-forming operations. The free structured category $\cF(\mathbb{G})$\/ generated by
$\mathbb{G}$\/ may then be taken to be $C(\cL_{\mathbb{G}})$, the syntactic category of $\cL_{\mathbb{G}}$.
We thus have the following picture:
\[
\xymatrix{
\mbox{\bf S-Cat} \ar@/_1ex/[rr]^\simeq_L  && \mbox{\bf Lang}  \ar@/_1ex/[ll]_C \\
& \mbox{\bf DirGrph} \ar@/_2ex/[ur]_{\mathcal{L}} \ar@/^2ex/[ul]^{\mathcal{F}}
}
\]
Of particular importance is the case where $\mathbb{G}$\/ is the empty graph. The resulting category
$\mathcal{F}(\mathbb{G})$\/ is then the initial structured category.
\end{enumerate}

In the book~\cite{LS86}, such theories include typed (and even untyped) lambda calculi 
(corresponding to cartesian closed categories with additional structure) and intuitionistic higher order logics 
(Russellian type theories) with full impredicative comprehension, extensionality, and Peano's axioms 
(corresponding to elementary toposes with logical morphisms and natural numbers).    We briefly discuss the two
cases of Cartesian Closed Categories and Elementary Toposes (both with natural numbers object) below. 

It is important to note that in order to obtain a 1-categorical equivalence $\mbox{\bf S-Cat} \simeq \mbox{\bf Lang}$\/ 
of this kind, we need to consider the objects of $\mbox{\bf S-Cat}$\/ not just as structured categories, but as 
categories \emph{equipped with specified structure}. Similarly, we require the functors to preserve this specified structure
on the nose. It is possible to avoid working with chosen structure, but then one should instead consider
$\mbox{\bf S-Cat}$\/ as a 2-category, and set up a 2-categorical equivalence with a suitable 2-category of theories.
An example of this finer analysis appears (in this volume) in the paper of Castellan et al.~\cite{CastellanS:catf}, which
discusses the Seely correspondence between locally cartesian closed categories and dependent type theories, and,
more generally, provides a suitable 2-categorical perspective on categorical logic.

\subsection{CCCs and the Curry-Howard-Lambek correspondence}\label{subsec:CH}

Cartesian closed categories were introduced by Lawvere in the early 1960s as the categorical analog of Church's 
typed lambda calculi.  In the early 1970s, Lambek explored this correspondence, along with connections to 
Sch\"onfinkel and Curry's works on combinatory algebras and functional completeness. 
 The precise tripartite categorical equivalence of cartesian closed categories, typed lambda calculi, and 
labelled deductive systems for positive intuitionistic propositional calculus (modulo equality of proofs)
was developed in detail in~\cite{LS86}.  This yields a modern version of the so-called Curry-Howard 
correspondence~\cite{GLT89}, with the additional idea (Lambek \cite{Lam68b,Lam69}) of equations between proofs, and is summarized in 
 Theorem~\ref{ccclam} below.

 \begin{definition} 
 A {\em cartesian closed category $\cC$ } (with specified structure) is a
 cartesian category $\cC$ (i.e., a category with specified finite products) such that, for each object $A\in
 \cC$, the functor
 $(-)\times A: \cC \rightarrow \cC$ has a specified right adjoint,
 denoted $(-)^A$. Thus, there is a natural  isomorphism (natural
 in $B \, \mbox{and} \, C$):
 \[ Hom_{\cC}(C\times A, B)   \cong  
 Hom_{\cC}(C, B^A).\]
 \end{definition}

\begin{example}{Examples of CCCs}
The category of sets is a CCC with $B^A$\/ the set of all functions $A \to B$. More generally, any functor category
$[\cC^{\mathrm{op}}, \set]$\/ is a CCC, where $G^F(C)$\/ is the set of natural transformations from 
$Hom_{\cC}(-,C) \times F$\/ to $G$. The category $\cat$\/ of small categories is also cartesian closed, 
as are many categories
of ``nice" topological spaces, such as compactly generated Hausdorff spaces. 
\end{example}

Next, consider simply typed lambda calculi. 
\begin{definition} A \emph{simply typed lambda calculus} consists of the following data.
First, it has a collection of \emph{simple types}
generated from a set of ground types $G$\/ by the grammar 
\[ {\bf Types}  \qquad A,B::= G \mid  \synt{1} \mid  A\times B \mid A\Rightarrow B.\]
At each type, we assume given an infinite set of variables; we write $x:A$\/ to indicate that
$x$\/ is a variable of type $A$.
Next, we have, for all types $A_1, \ldots, A_k,B$\/ a (possibly empty) set of basic terms $E(A_1, \ldots, A_k; B)$.
Then the collection of typed \emph{terms} is generated 
using the rules displayed in Figure~\ref{fig:term}. 
We make the usual assumptions (see e.g.~\cite{BarendregtH:lamcss, LS86})
regarding free and bound variables, and write $FV(t)$\/ for the set of 
free variables of $t$; each $x \in FV(t)$\/ has a unique type, and from the term $t$\/ we
can recover the types of the free variables in $t$.
\begin{figure}[h]
$$
\begin{tabular}{|p{5pt}cccp{5pt}|}
\hline 
&& & &\\[-5pt]
&$\infer{*:\synt{1}}{}$ & \ \ $\infer{f(x_1, \ldots, x_k):B}{f\in E(A_1, \ldots, A_k;B) & x_i: A_i}$ \ \   
&\ \ $\infer{\lambda x:A.\varphi(x): B^A }{x: A & \varphi(x): B}$
& \\[2ex] 
&$\infer{\langle a, b\rangle: A\times B}{a: A & b: B}$& $\infer{\pi_it: A_i}{t: A_1 \times A_2}$&   $\infer{ft:B}{f:B^A & t: A}$ & \\[1ex]   
&& && \\
\hline
\end{tabular}
$$
\caption{Lambda Calculus Terms}
\label{fig:term}
\end{figure}

\noindent
 Finally, we have equations between terms of the same type. 
We write \\$t=_Xs:A$\/ to express that the terms $t,s$\/ are equal, and 
that the free variables of $t$\/ and $s$\/ are contained in the set $X$.

The relations $=_X$\/ are congruences
satisfying the following clauses \footnote{Here we present lambda calculi as ordinary equational theories, as in
 \cite{LS86}. One could also write equational logics in an appropriate sequent calculus, writing
 $t=_Xs:A$  as $\vdash_X s = t:A$ (cf. Barendregt's {\em lambda theories} \cite{BarendregtH:lamcss} 
 and the use of HOL below).}:
\begin{itemize}
\item $t=_Xs$, $X \subseteq Y$\/ implies $t=_Ys$
\item $t=_Xs$\/ implies $ft =_X fs$ (where $f:B^A$\/ and $t,s:B$)
\item $\varphi(x)=_{X \cup \{x\}} \psi(x)$\/ implies $\lambda x.\varphi(x) =_X \lambda x.\psi(x)$ 
\item $a=_X*$\/ (where $a:\synt{1}$)
\item  $\pi_i\langle a_1,a_2\rangle=_Xa_i \; ; \quad a=_X\langle \pi_1a,\pi_2a\rangle$
\item $(\mathbf{\beta})$ \quad $(\lambda x.\varphi(x))t=_X\varphi[t/x] \; ; \qquad (\eta) \quad 
f=_X\lambda x.fx \; \mbox{ where } x \not \in FV(f)$
\end{itemize}
\end{definition}

 It is possible to augment simply typed lambda calculus with additional types, terms, and equations
(cf. \cite{LS86}).  We discuss the case of adding natural numbers and lists in Subsection \ref{subsec:NNOs} below. 


An important example of a simply typed lambda calculus arises as follows.

\begin{definition}[Simply typed $\lambda$-calculus from a graph]
Consider a directed graph $\mathbb{G}=(G,E)$.
The calculus $\mathcal{L}_\mathbb{G}$\/ has as ground types the vertices of $\mathbb{G}$,
and as basic terms the edges of $\mathbb{G}$,  (i.e. whenever $f:A \to B$\/ is in $E$, there is a basic term
$f(x):B$, with $x:A$.) The congruence $t=_Xs$\/ on terms is the smallest
congruence satisfying the rules of simply typed lambda calculus.
\end{definition}

We now define the category ${\bf CCC}$\/ whose objects are CCCs (with chosen products and exponentials),
and whose morphisms are functors preserving the chosen products and exponentials on the nose.
On the other hand, we define the category {\bf Typed $\lambda$-calc} to have typed lambda calculi as
objects, and translations as morphisms. Here, a \emph{translation} between two calculi is a mapping sending types to 
types and terms to terms, in such a way that all type and term formation operations are preserved and that
provable equality between terms is preserved. 

\begin{definition}[Internal language of a CCC]
Let $\cC$\/ be a cartesian closed category.  The \emph{internal language} of $\cC$\/ is the simply typed lambda 
calculus $L(\cC)$\/ generated by the underlying graph of $\cC$, together with all equations holding between arrows of $\cC$.
\end{definition}

In the other direction, we construct a CCC $C(\cL)$\/ from a typed lambda calculus $\cL$:

\begin{definition}[Syntactic Category]
\label{syntactic_ccc}
Let $\cL$\/ be a simply typed lambda calculus. Define a category $C(\cL)$\/ by:
\begin{description}
\item[{\bf Objects}] The types of $\cL$.
\item[{\bf Morphisms}] For any term $t:T$\/ with $FV(t)=\{x_1:T_1, \ldots, x_n:T_n\}$, we have a morphism
$[t]:T_1 \times \cdots \times T_n \to  T$. Here $[t]$\/ is the equivalence class of $t$\/ under provable
equality of the theory $\cL$.
\item[{\bf Identities}] The identity at an object $T$\/ is represented by the term $x:T$.
\item[{\bf Composition}] Given terms $t(x),s(y)$\/ representing morphisms $A \to B$\/ and $B \to C$\/ respectively
(where we assume that $t$\/ is substitutable for $y$\/ in $s$),
the term $s[t/y]$\/ (the result of substituting $t(x)$ for all variables $y$ \/ in $s$) represents the composite $A \to C$. 
\end{description}
\end{definition}
We now have the promised result\footnote{Lambek reported that
when he lectured at Columbia announcing these results  
Sammy Eilenberg is reported to have said:  ``This is wonderful.  
 Now category theorists will never have to learn lambda calculus!" }:

 \begin{theorem}[Curry-Howard-Lambek correspondence~\cite{LS86}]\label{ccclam}  The pair of functors $L:\mbox{\bf CCC} \to 
 \mbox{\bf Typed $\lambda$-calc}$\/ (internal language) and
 $C: \mbox{\bf Typed $\lambda$-calc} \to \mbox{\bf CCC}$\/ constitute an equivalence of categories.
  \end{theorem}
  
 The above theorem extends to include adding the natural numbers and 
 similar data types (of which the categorical aspects are discussed in the next Section). 

\subsection{Elementary toposes and HAH}\label{subsec:HAH}
We now outline another instance of an equivalence between a class of categories and
a fragment of logic, namely elementary toposes with NNO and higher-order intuitionistic arithmetic (HAH).

Recall that in a category $\cC$, a \emph{subobject} of an object $A$\/ is an equivalence class of monomorphisms
$m:X \to A$, where two monomorphisms are equivalent precisely when they factor through each other. 
The collection of subobjects of $A$\/ is denoted $Sub(A)$. The assignment $A \mapsto Sub(A)$\/ is a contravariant 
functor from $\cC$\/ to the category of posets.

A category is said to have \emph{canonical} subobjects when every subobject has a chosen representative. In $\set$,
for example, we may represent a subobject through its set-theoretic image. 
 
\begin{definition}[Elementary Topos]
A category $\cC$\/ is a \emph{topos} when it has the following structure:
\begin{itemize}
\item $\cC$\/ has finite limits
\item $\cC$\/ has power-objects: for each $A$\/ there exists an object $\synt{P}A$\/ and natural bijection
\[ Hom_{\cC}(B,\synt{P}A) \cong Sub(B \times A).\]
\end{itemize}
\end{definition}
A power-object $\synt{P}A$ (when it exists) represents the functor $Sub(- \times A)$. In the category of sets,
we may take $\synt{P}A$\/ to be the powerset of $A$, and then the defining bijection becomes the familiar correspondence between
relations $R \subseteq B \times A$\/ and functions $r:B \to \synt{P}A$.

In a topos, we write $\synt{\Omega}$\/ for $\synt{P}1$. This is the \emph{subobject classifier}: there is a natural bijection
\[ Hom_{\cC}(B,\synt{\Omega}) \cong Sub(B). \]
We think of $\synt{\Omega}$\/ as the object of truth values of $\cC$, and of $\synt{P}A$\/ as the exponential $\synt{\Omega}^A$.

The category of sets is of course a typical example of a topos, 
as are functor categories $[\cC^{\mathrm{op}},\set]$. Other
examples will be discussed below. The qualifier \emph{elementary} is used to stress the inclusion of toposes other
than \emph{Grothendieck Toposes} (which are required to be cocomplete and have a small set of generators). 
\footnote{Grothendieck toposes were introduced in the early 1960s 
by the Grothendieck school of algebraic geometry~\cite{SGA} as sheaves on 
a site. In the early 1970s, Lawvere and Tierney~\cite{ICM70} introduced elementary toposes. It was realized that such toposes
could be considered as a universe of mathematics, where the objects and morphisms can be treated as sets and functions,
provided one refrains from using classical reasoning (the law of excluded middle and the Axiom of Choice).}
 
In the context of elementary toposes, one often considers \emph{logical morphisms} between toposes. These
are functors preserving all the topos structure. Just as for CCCs, we work with toposes with specified structure and 
morphisms strictly preserving this structure. 

\begin{definition}
The category ${\bf Top}$\/ has:
\begin{description}
\item[{\bf Objects}] Elementary toposes with specified finite limits and power objects, and with canonical subobjects.
\item[{\bf Morphisms}] Logical functors preserving all specified structure on the nose.
\end{description}
\end{definition}

Next, let us describe (intuitionistic) \emph{higher-order logic} (HOL). 
This formal system can be thought of as an extension of simply typed lambda calculus, 
with added type and term constructors for the type $\synt{\Omega}$\/ of propositions and for power objects $\synt{P}A$. 
(However, we do not include exponentials explicitly, as they are definable in terms of the other operations
\footnote{Moreover, as discussed in \cite{LS86}, strict logical functors will preserve only the powerset structure on the nose.
In keeping with the logic literature and because of its historical importance, we denote the type of truth values  
by $\synt{\Omega}$, rather than treating it as $\synt{P}1$. Logical functors will preserve $\synt{\Omega}$ on the nose. }.)
Thus the types are generated from ground types $G$\/ using the following grammar:
\[ {\bf Types} \qquad A,B ::= G \mid \synt{1} \mid A \times B  \mid \synt{\Omega} \mid \synt{P}A. \]
The terms are generated from basic terms and variables using the rules displayed in Figure~\ref{fig:HAH}
(where we omit the rules already stated for simply typed lambda calculus in Figure~\ref{fig:term}):
\begin{figure}[ht]
\[
\begin{tabular}{|ccccc|}
\hline 
& \infer{a=a' : \synt{\Omega}}{a:T & a':T} \ \ & \ \ \infer{a \in \alpha : \synt{\Omega}}{a:T & \alpha: \synt{P}T} 
& \ \ \infer{\{x:A \mid \varphi(x)\}:\synt{P}A}{x:A  & \varphi(x): \synt{\Omega}} & \\
\hline
\end{tabular}
\]
\caption{Terms of higher-order intuitionist logic}
\label{fig:HAH}
\end{figure}

In~\cite{LS86} there are two axiomatizations of higher order logic, including the one above based on 
 {\em equality (between terms of the same type), comprehension, extensionality, and (in case 
 we add a type of natural numbers) Peano's axioms}.  

 Following Russell, Henkin, and Prawitz, since we are assuming a primitive equality
 predicate  at each type, we can define the usual logical connectives as in Figure \ref{fig:type} below.
\begin{figure}[ht]
\begin{center}
$\begin{array}[t]{llllll}
\top & := \quad & * = * & p\vee q & := \quad & \forall_{x:\synt{\Omega}}(((p\Rightarrow x) \wedge (q\Rightarrow x))\Rightarrow x)\\
p\wedge q & := & \langle p,q\rangle = \langle \top, \top\rangle \qquad & \forall_{x:A}\varphi(x) & := & \{x : A \ | \ \varphi(x)\} = \{x:A \ | \ \top\} \\
p\Rightarrow q & := & p\wedge q = p          & \exists_{x:A}\varphi(x) & := & \forall_{y:\synt{\Omega}}(\forall_{x:A}((\varphi(x)\Rightarrow y)\Rightarrow y))\\
\perp              & : = & \forall_{x:\synt{\Omega}} x                        & \exists !_{x:A}\varphi(x) & := & \exists x':A (\{x:A \ | \ \varphi(x)\} = \{x:A \ | \ x = x' \} )\\
\neg p             & := & \forall_{x:\synt{\Omega}} (p\Rightarrow x) &  & &                    \\
\end{array}$
\caption{Type-theoretic encoding of logic}
\label{fig:type}
\end{center}
\end{figure}

We now define an entailment relation $\Gamma \vdash_X q$. Here, $\Gamma$\/ is
a finite set of formulas (i.e., terms of type $\synt{\Omega}$), 
$q$\/ is a formula, and $X$\/ is a typed variable context containing all the
free variables of $\Gamma$\/ and $q$; the meaning of $\Gamma \vdash_X q$\/ is that
$q$\/ can be derived (using the rules for intuitionistic logic) from $\Gamma$.  When $\Gamma=\emptyset$\/ 
we simply write $\vdash_X q$. There are standard
structural rules (including Cut), substitution, rules for equality, rules for products, and for comprehension. For example, 
there is the comprehension rule
\[ \vdash_X (y \in \{x:A \mid \varphi(x) \}) = \varphi(y). \]
We refer to~\cite{LS86} for a complete list of rules. 

By a \emph{type theory} we mean an extension of HOL by sequents $\Gamma \vdash_X q$. When $\mathcal{L}$
is such a type theory, we write $\vdash^{\mathcal{L}}$\/ for the entailment relation of $\mathcal{L}$\/ (although we may 
omit the superscript when $\mathcal{L}$\/ is understood). In $\mathcal{L}$, we say that two terms $t,s$\/ of the same type
are provably equal when $\vdash^{\mathcal{L}}_X  t=s$. Just as for simply typed lambda calculi, it is common 
to include a type of natural numbers; the type theory obtained by adding the natural numbers to HOL (and no
further basic types) is called \emph{Higher-order intuitionistic Arithmetic}, or HAH for short. 

An \emph{interpretation} of one type theory in another is a mapping of types to types that preserves all type
formation operations, together with a mapping of terms that respects the typing, the term formation operations and
the provable equality. Type theories and interpretations form a category denoted $\mathbf{Lang}$.

A type theory is  {\em classical}
 if in addition it has Aristotle's axiom of excluded middle:  $\forall p: \synt{\Omega} (p \vee \neg p).$ Such a system
 of classical type theory was employed in G\"odel's famous incompleteness paper \cite{GodelK:ubefus}.

Given a type theory $\cL$ one may now build a syntactic topos as follows:
\begin{definition}[Generated Toposes $T(\cL)$] 
\label{topgen}  The  topos {\em $T(\cL)$ generated by the type theory $\cL$} has as objects ``sets" (i.e., closed terms
$\alpha$ of type $\synt{P}A$ , modulo provable equality). Morphisms
$\alpha\rightarrow \beta$, where $\alpha:\synt{P}A$ and $\beta:\synt{P}B$,
are   ``provably functional relations", i.e. closed terms $\varphi:\synt{P}(A\times B)$
(modulo provable equality) such that we can prove:
\[
\vdash^{\cL} \forall_{x: A}(x\in \alpha \Rightarrow \exists!_{y: B} (y\in\beta \wedge (x,y)\in
\varphi))
\]
$T(\cL)$ is the category of ``sets"  and ``functions" 
formally definable  in higher-order logic $\cL$. 
\end{definition}

The assignment $\cL \mapsto T(\cL)$\/ is a functor $T:\mathbf{Lang} \to \mathbf{Top}$.
For $\cL_0$ = pure type theory, $T(\cL_0)$ is called the {\em free topos}, denoted $\Ftop$. It enjoys the following 
universal property: for any elementary topos $\cE$\/ there exists a logical functor $F:\Ftop \to \cE$\/ which is unique
up to (unique) natural isomorphism. In other words, $\Ftop$\/ is the initial object of $\mbox{\bf Top}$. 

In the other direction we may assign to a topos $\cC$\/ its \emph{internal language} $L(\cC)$, just as for
CCCs. This gives a functor $L:\mathbf{Top} \to \mathbf{Lang}$.

\begin{theorem}[Lambek-Scott~\cite{LS86}]
The functors $L,T$\/ described above constitute an equivalence of categories
\[
\xymatrix{
\mbox{\bf Top} \ar@/_1ex/[rr]_L^\simeq && \mbox{\bf Lang}. \ar@/_1ex/[ll]_T \\
}
\]
\end{theorem}

As for simply typed lambda calculus, we may extend this result by adding datatypes. Most importantly,
we can consider type theories with natural numbers and toposes with natural number objects (see next Section).

\section{What are computable functions in categories?}\label{sec:compcat}
We turn to the study of computable functions in categories. In this section, we limit ourselves to 
computable numerical functions; later we shall consider computable maps on other datatypes. 

\subsection{Natural Numbers Objects and \prim}
\label{subsec:NNOs}
 In order to discuss number-theoretic functions in categories, we  briefly recall Lawvere's notion of 
 {\em Natural Numbers Objects} (NNOs) in cartesian closed categories~\cite{Law64, LS86} and more 
 generally NNOs in cartesian and monoidal categories~\cite{PareR:moncnn}.
 
 \begin{definition}[Lawvere~\cite{Law64}] A {\em Natural Numbers Object} (NNO) in a (cartesian closed) category  $\cC$
 is a diagram $\synt{1}\stt{0}\csynt{N}\stt{S}\csynt{N}$ initial among diagrams 
 $\synt{1}\stt{a}A\stt{h}A$.  
 i.e., there exists a unique $It_{ah}: \csynt{N}\rightarrow A$
 satisfying:    $$It_{ah} \circ 0 = a \ \ , \ \ It_{ah} \circ S = h \circ It_{ah}$$ 
 Existence, without uniqueness,
 of such an arrow $It_{ah}$ yields the notion of a {\em weak NNO}.  Any arrow $It_{ah}:  \csynt{N}\rightarrow A$ (unique or not) satisfying the
 equations above is called
an {\em iterator at type $A$}. Diagrammatically, 
 
  \end{definition}

\qquad  

\begin{minipage}[b]{2in}
 \adjustbox{scale=1.1}{%
    \begin{tikzcd}
    \synt{1}\arrow{r}{0}\arrow{dr}{a} & \csynt{N} \arrow[d,dotted,"It_{ah}"]\arrow{r}{S}& \csynt{N} \arrow[d,dotted,"It_{ah}"] \\
  & A\arrow{r}{h}& A\\
\end{tikzcd}
}\end{minipage}
\qquad \begin{minipage}[b]{3in}
In $\set$ this says:\\[1ex] 
$It_{ah} (0)  =  a  \\
  It_{ah}(n+1) =  h(It_{ah}(n))$
\end{minipage}

 
 
For any NNO (weak or strong) we may define, for any natural number 
$n \in \mathbb{N}$, the \emph{standard numeral} $\numer{n}:\synt{1} \to \csynt{N}$\/ by
\[ \numer{0}=0 \; \quad \numer{n+1}=S \circ \numer{n}.\]
We stress that depending on the nature of the ambient category, there may be non-standard numerals, that is, points $\synt{1} \to \csynt{N}$\/ 
that are not of the form $\numer{n}$.

A natural numbers object in a cartesian closed category is equivalent to the following
scheme of {\em Iteration with parameters}.  This general
scheme (and its variants for monoidal categories) is sufficient for representing the 
primitive recursive functions~\cite{Law64,LS86} and is the appropriate definition for NNO's
in cartesian (or monoidal) categories, as in~\cite{LS86}, p.71.

\begin{definition}[Parametrized NNO]
\label{paramnno}
  A diagram $\synt{1}\stt{0}\csynt{N}\stt{S}\csynt{N}$ in a cartesian
category $\cC$ is a {\em parametrized NNO} if for all arrows $A\stt{g} B, B\stt{f}B$,  there exists a unique $It_{gf}: \csynt{N}\times A\rightarrow B$ such that:  
$$\begin{minipage}[b]{2in}
 \adjustbox{scale=1.1}{%
    \begin{tikzcd}
   A\arrow{r}{\langle 0!, 1_A\rangle}\arrow{dr}{g} & \csynt{N}\times A \arrow[d,dotted,"It_{gf}"]\arrow{r}{S\times 1_A}& \csynt{N}\times A \arrow[d,dotted,"It_{gf}"] \\
  & B\arrow{r}{f}& B\\
\end{tikzcd}
}
\end{minipage} \qquad \qquad \qquad 
\begin{minipage}[b]{2in}
In $\set$ this says: \\[1ex]
$It_{gf} (0,a)  =  g(a) \\  
It_{gf}(n+1,a)  =  f(It_{gf}(n,a))$
\end{minipage}  
$$

\vspace{-3ex}

\noindent
Existence without uniqueness of the arrow $It_{gf}$ above yields a {\em weak parametrized NNO}.
\end{definition}

A typical example is the notion of adding an iterator to a simply typed lambda calculus.

\begin{example}{Iterators in typed lambda calculus}
Following~\cite{LS86}, we add to the terms of simply typed lambda
calculus in Figure \ref{fig:term} an atomic type $\synt{N}$ and term
formation operations
\[ \synt{0}:\synt{N} \qquad \infer{\synt{S}x:\synt{N}}{x:\synt{N}} \qquad \infer{\iter(a,h,x):A}{a:A & h:A^A & x:\synt{N}} \]
(allowing in particular the definition of standard numerals $\numer{n}$).
We then add to the equations of the simply typed lambda calculus the following 
equations:
\[
\iter(a,h,0) =_X a    \qquad   \iter(a,h,\synt{S}x) =_{X \cup \{x\}} h(\iter(a,h,x)), 
\mbox{provided $x\not\in X$}.
\]
Calling this lambda theory $\cL$, the associated syntactic category $C(\cL)$ (Definition \ref{syntactic_ccc}) 
is a cartesian closed category with weak NNO.

\end{example}

In general, when we consider a category $\cC$\/ 
the difference between a weak and a strong NNO in $\cC$\/ can be understood in logical terms by considering the form of induction
allowed in the internal language. For example, when $\cC$\/ has a strong NNO we can prove the entailment
\[ \vdash_{x,y} x+y=y+x \]
where, crucially $x,y$\/ are free variables of type $\csynt{N}$. When $\cC$\/ only has a weak NNO one can prove by 
(external) induction that for every $n \in \mathbb{N}$:
\[ \vdash_{x} x+\numer{n}=\numer{n}+x. \]

Next, consider a (not necessarily symmetric) 
monoidal category $(\cC, \otimes, I)$.  
Following Par\'e and Rom\'an~\cite{PareR:moncnn}, we may define notions of 
{\em Left} and {\em Right} NNOs, in analogy with Definition \ref{paramnno}.

\begin{definition}[Left Parametrized NNO] A diagram $I \stt{0}\csynt{N}\stt{S}\csynt{N}$ in 
a monoidal category $\cC$ is a {\em left parametrized NNO} if for all arrows  $A\stt{g} B, B\stt{f}B$,  
there exists a unique $k: \csynt{N}\otimes A\rightarrow B$ such that:  
$$
 \begin{minipage}[b]{4in}

 \adjustbox{scale=1.1}{%
    \begin{tikzcd}
 I\otimes  A\arrow{r}{0\otimes 1_A}\arrow{d}{\cong} & \csynt{N}\otimes A \arrow[d,dotted,"k"]\arrow{r}{S\otimes 1_A}& \csynt{N}\otimes A \arrow[d,dotted,"k"] \\
 A\arrow{r}{g} & B\arrow{r}{f}& B\\
\end{tikzcd}
}

\end{minipage}
$$
\end{definition}

\vspace{-2ex}

\noindent
In the same manner, tensoring by $A$ on the left (rather than the right) results in a {\em Right
Parametrized NNO};  {\em weak} objects are defined similarly by  assuming merely 
existence (but not necessarily uniqueness) of $k$.  For many examples of such monoidal NNOs, see~\cite{PareR:moncnn}.

We remark that there are yet other axiomatizations. A \emph{Peano-Lawvere} category
is a category for which the forgetful functor $\cC^\NN \to \cC$\/ has a left adjoint (where $\NN$\/ is regarded as the 
free monoid on one generator). A systematic study of the free such category can be found in Burroni's ~\cite{BurroniA:recgI}.

Another relevant class of categories is that of \emph{list-arithmetic pretoposes}. These were developed by
Maietti~\cite{MaiettiM:joyaul} (see also~\cite{MaiettiM:indpca})
 in order to provide a categorical setting accommodating Joyal's \emph{arithmetic universes} (\cite{Joyal05}),
which in turn serve as a categorical account of the Incompleteness Theorem. A pretopos is a category that has finite limits,
pullback-stable disjoint coproducts, and pullback-stable quotients of equivalence relations. Such a category has
\emph{parameterized list objects} when for each object $A$\/ there is an object $LA$\/ equipped with maps
$e:\synt{1} \to LA, c:LA \times A \to LA$\/ (thought of as the empty list and concatenation). These are required to satisfy the
following universal property: for any $a:B \to C$\/ and $h:C \times A \to C$\/ there is a unique $It_{a h}:B \times LA \to C$\/
 making the following commute:
\[
\xymatrix{
B \ar[dr]_a \ar[r]^-{\langle 1_B,e\rangle} & B \times LA \ar@{-->}[d]^{It_{a h}} && B \times LA \times A \ar[ll]_{1_B \times c} 
\ar@{-->}[d]^{It_{a h} \times 1_A}  \\
& C && C \times A \ar[ll]_h
}
\]
As is the case for NNOs, we may also consider a weak version where we only require existence and not uniqueness
of the iterator $It_{ah}$.
A list-arithmetic pretopos is now defined as a pretopos admitting parametrized list objects for all $A$. Note that taking $A=\synt{1}$\/ 
gives the notion of a parameterized NNO.

\subsection{Representability}\label{subsec:repr}
We now turn to representability of numerical functions in categories with NNOs. 

\begin{definition}[Lambek-Scott~\cite{LS86}]  
\label{rep}
Let $\cC$ be a cartesian category with a weak 
parametrized NNO 
 \ $\synt{1}\stt{0}\csynt{N}\stt{S}\csynt{N}$.  A function
$f:\NN^k\rightarrow \NN$ is {\em representable} in $\cC$ if there is an arrow
$F: \csynt{N}^k\rightarrow \csynt{N} \in \cC$ such that 
$F\langle \numer{n_1}, \cdots , \numer{n_k}\rangle = \numer{m}$
whenever $f(n_1,\cdots,n_k)=m$. 
 \end{definition}
 
 \noindent

 Of course, the determination of which numerical functions are representable depends
 on the category:  in the category  $\set$, {\em all} numerical functions are representable!
 Following Lambek's question in the Introduction, we shall look at free categories with NNOs.
 
 
  \begin{theorem}[Rom\'an~\cite{RomanL:carcnn}] 
 The class of representable numerical functions in $\Fcart$, the free cartesian category  with 
 parametrized NNO,   is $\prim$.
 \end{theorem}
 
 \noindent
 Hence the unique representation functor $\Fcart\longrightarrow \set$ has
 as image the subcategory of sets whose objects are powers $\NN^n$ and whose maps are
 tuples of primitive recursive functions. 

Rom\'an's proof essentially shows that Goodstein's 
development~\cite{GoodsteinR:RNT} of Skolem's primitive recursive arithmetic can be mimicked
in $\Fcart$.  In that sense, the result is not so surprising.   However the following striking result considers
the extension to $\Fmon$, the free monoidal category with a LNNO. 
Recall primitive recursion requires projection functions $U^n_i: \NN^n\rightarrow \NN$, 
yet in a monoidal category, in general  $\otimes$ does not have explicit projections. Nevertheless:
  \begin{theorem}[Par\'e-Rom\'an~\cite{PareR:moncnn}] \
  \begin{enumerate}[(i)]
\item The primitive recursive functions are representable in any monoidal category with LNNO.

\item Indeed, $\Fmon$, the free monoidal category with LNNO, exists and is isomorphic to $\Fcart$, 
the free cartesian category with parametrized NNO.
\end{enumerate}
\end{theorem}

\noindent
Why is this?  The reason is that the objects of $\Fmon$ are generated from $\csynt{N}$ under tensoring and  
{\em we can code projections and diagonals between tensor powers $\csynt{N}^{\otimes k}$}. 
This then allows the representability of the primitive recursive functions  in a similar manner to 
$\Fcart$.  The former result (coding projections and diagonals) is proved by an elegant categorical 
argument in Par\'e-Rom\'an,  while Plotkin~\cite{Plotkin13} gives a direct (albeit nontrivial) coding of the primitive 
recursive functions in $\Fmon$.

\subsection{Going beyond the primitive recursive functions: free CCCs}
\label{subsec: freeCCCs}
How do we get more representable functions?  We increase the logical strength (the types)
from the logic of $\{\wedge, \top\}$ (or $\{\otimes , I\}$) to the logic of $\{\wedge, \Rightarrow, \top\}$, i.e. to 
the cartesian closed level.  Consider the free CCC with natural numbers generated by
the empty graph, denoted $\Fccc$\/ (as defined in Section~\ref{subsec:CH}).

The following is a theorem about simply typed lambda calculus, translated into
the language of CCCs:  
\begin{theorem}[Lambek-Scott~\cite{LS86}]
In $\Fccc$, the free CCC with weak NNO:
\begin{enumerate}
\item All primitive recursive functions and the Ackermann function are representable.
\item The representable functions form a proper subclass of the total recursive functions, namely
the provably total functions of Peano Arithmetic, or equivalently the $\varepsilon_0$-recursive 
functions~\cite{TroelstraA:baspt,GLT89}.
\end{enumerate}
\end{theorem}

\noindent
More generally, the representable total functions of $\Fccc$\/ are the lowest level of the hierarchy
of G\"odel's {\em Dialectica Functionals}, i.e., G\"odel's primitive recursive functionals
of finite type~\cite{GodelK:ubeebu,TroelstraA:metiia}.

There is also a version of G\"odel's Incompleteness for $\Fccc$\/.
 Let $Z$  represent the zero function.

\begin{theorem}[A version of Incompleteness, or $\synt{1}$ is not a generator] In $\Fccc$, there is a 
  closed term $F: \synt{N}\Rightarrow \synt{N}$ such that for each numeral $\numer{n}$,
 $\vdash F\numer{n} = \numer{0}$, but \,
 $\not \vdash F = Z$ .
  \end{theorem}
For a proof, see  Corollaries 2.11, 2.12 in \cite{LS86}, p.263.

Finally, a topic of considerable importance in theoretical computer science:
\begin{example}{Computation by normalization}
We should also recall the notion of {\em computation by normalization or, for a logician, by cut-elimination} 
~\cite{GLT89}.  In the rewriting theory of typed
lambda calculus, we can set up strongly normalizing rewrite systems in which
terms can be rewritten to (unique) normal forms.  
 
Given a term $f: \synt{N}\Rightarrow \synt{N}$ and a numeral $\numer{n}:\synt{N}$, to 
compute $f\numer{n}$ by normalization,  we first normalize
this term to its unique normal form of type $\csynt{N}$.  This yields a numeral $\numer{m}$, for which we can prove 
$\vdash f{\numer{n}} = \numer{m}$; cf. \cite{GLT89}. This gives the value of $f$ on input numerals.

By Curry-Howard-Lambek, normalization techniques may also be used to solve coherence problems 
(decidability of equality) for various free CCCs, via their internal languages \cite{LS86}: to check if 
two arrows in a free CCC are equal or not, it suffices to show that their  normal forms (qua lambda terms) are identical,
up to change of bound variables.    
  
 \end{example}

 \subsection{Some properties of the free topos}   

 Pure intuitionistic type theory $\cL_0$  has many interesting properties, which translate into
 algebraic properties of the free topos $\Ftop$ (see~\cite{LS86}) and are often key metamathematical 
 principles of intuitionistic systems (\cite{TroelstraA:metiia}). In what follows we write $\vdash$\/ instead
 of $\vdash^{\cL_0}$\/ for derivability in intuitionist higher order arithmetic HAH.

\begin{itemize}\itemsep=1ex
\item {\bf Consistency:}  not ( $ \vdash \perp$) .
\item {\bf Disjunction Property:} If \, $\vdash p\vee q$ , then
$\vdash p$ or
$\vdash q$.
\item {\bf Existence Property (EP):}  If $\vdash \exists_{x:A}\varphi(x)$ then 
$\vdash \varphi(a)$ for some closed term $a : A$.

\noindent
In particular, in $\Ftop$  EP says that numerals are standard, i.e. that
numerals $\synt{1}\stt{f}\csynt{N}$ are all of the form $\numer{n}$, for some $n\in\NN$.

\item {\bf Troelstra's Uniformity Principle (UP) for $A = \csynt{P}C$}:   \\[1ex]
If \ $\vdash \forall_{x:A}\exists_{y:\csynt{N}}\varphi(x,y)$ then 
$\vdash  \exists_{y:\csynt{N}}\forall_{x:A}\varphi (x,y)$.

\noindent
  In $\Ftop$,   UP says the arrows $\csynt{P}C
\rightarrow \csynt{N}$ are constant, i.e. factor through a standard numeral.

\item {\bf Independence of premisses (IP): } If \ $\vdash \neg p \Rightarrow
\exists_{x:A}\varphi(x)$ then \\ $\vdash \exists_{x:A}(\neg p\Rightarrow  \varphi(x))$.


\item {\bf Markov's Rule  (MR): } If \ $\vdash \forall_{x:A}(\varphi(x)\vee \neg \varphi(x))$
and $\vdash \neg\forall_{x:A}\neg \varphi(x)$, then  $\vdash \exists_{x:A}\varphi(x)$.

 \item  {\bf The Existence Property with a parameter of type $A =  \csynt{P}C$: } \\
If  \ $\vdash \forall_{x:A}\exists_{y:B}\varphi(x,y)$ then $\vdash \forall_{x:A}\varphi(x,\psi(x))$, 
where $\psi(x):B$.

\end{itemize}

\noindent
{\bf Proofs:} The original proofs~\cite{LS80}  for EP and DP used Friedman (impredicative) realizability.  When the authors
lectured on this, Peter Freyd realized these rules had purely algebraic statements, with direct categorical proofs, using Artin gluing categories
(\cite{Wraith74}).  The Freyd gluing techniques were expanded to include the proof rules above in \cite{LS80} and in a series of later papers
by the authors.  This is also presented in \cite{LS86}.


The {\em free Boolean topos} $\Fbtop$ is defined in the same way as the free topos, but generated from 
{\em classical} type theory.  As argued in~\cite{LS86}, alas the free Boolean topos is {\em not} an ideal universe 
for classical mathematicians. For example, as a consequence of G\"odel's Incompleteness Theorem, 
there are non-standard numerals.  To see this, let $G$ be 
any undecidable G\"odel sentence.  It may be shown that 
$\varphi(x):= (x = 0 \Rightarrow G) \, \wedge \, (x\not = 0 \Rightarrow \neg G)$ determines a numeral 
$f: {\synt{1}} \rightarrow \csynt{N}$ in $\Fbtop$; however, it cannot be a standard numeral, else we could decide $G$.
 
\vspace{1ex}
We now turn to the matter of representable numerical functions in the free topos.
First we recall the definition of representability of a function in HAH\footnote{The same definition works in 
other formal systems such as Peano Arithmetic.}:

\begin{definition}[Representability in HAH]
A total function $f:\NN^k \to \NN$\/ is \emph{representable} in HAH when there exists a formula 
$R_f(x_1, \ldots, x_k,y)$\/ 
such that
\begin{enumerate}[(i)]
\item $f(n_1, \ldots, n_k)=m$\/ if and only if $\vdash R_f(\numer{n_1}, \ldots, \numer{n_k}, \numer{m})$
\item $\vdash \forall x_1:\synt{N} \ldots x_k:\synt{N} \exists ! y:\synt{N}.R_f(x_1 \ldots, x_k,y)$.
\end{enumerate}
\end{definition}
In the literature, one often considers a weaker notion of representability, in which clause (i) remains,
but (ii) above is
replaced by 

\medskip
\noindent
 (ii')  for all  $n_1, \ldots, n_k \in \NN. \vdash \exists ! y:\synt{N}. R_f(\numer{n_1}, \ldots, \numer{n_k}, y)$.
 
\medskip
\noindent 
We refer to this weaker notion as \emph{numeralwise representability}.

It follows that a total numerical function  $\NN^k \to \NN$ is representable in HAH if and only if it is representable by an arrow
$N^k\to N$  in the free topos. 
(See Prop. 3.1, p. 264 in~\cite{LS86} for details.)

\begin{theorem}[Lambek-Scott~\cite{LS86}]  ~~
\begin{enumerate}[(i)]
\item In HAH (and hence in the free topos), every representable numerical function
 is recursive. In particular, the global sections functor $\Ftop(\synt{1},-): \Ftop \to \set$\/ sends
morphisms $\csynt{N}^k \to \csynt{N}$ to recursive functions $\NN^k \to \NN$. 
\item Not all total recursive functions so arise.
\end{enumerate}
\end{theorem}
(The second part of the theorem can be established by means of a diagonal argument.)
This of course leads to the question of which total recursive functions are representable in HAH. 
This is related to the representability of numerical functions in Girard's system $\mathcal{F}_\omega$, but we shall
not pursue this here.
We note that the situation changes radically if we consider {\em classical} type theory (the free Boolean topos).
 
\begin{theorem}[Lambek-Scott~\cite{LS86}] ~~
\begin{enumerate}[(i)]
\item The numeralwise representable functions in classical type theory are exactly the
total recursive functions  (G\"odel).
\item In classical type theory,   numeralwise representable functions coincide with representable ones
(by a result of V. Huber-Dyson \cite{VHD65}, detailed in ~\cite{LS86}).    Hence the representable
numerical functions in classical type theory are exactly the total recursive functions.
 
\end{enumerate}
\end{theorem}
Unfortunately, as we have seen, the free Boolean topos has non-standard numerals.  Thus, the global 
sections functor from the free Boolean topos to \set\/ in general sends arrows $\csynt{N}^k\to \csynt{N}$ to partial, 
rather than total, numerical functions.
This suggests that the representability of {\em partial functions} may be at least
as important as that of total functions. In fact, we shall see that even at the intuitionistic level, the theory becomes
much smoother.

\begin{definition}  
\label{parrep}
A partial function $f: \NN^k\rightharpoonup\NN$ is {\em representable} in  HAH if there is a formula 
$R_f(x_1, \ldots, x_k,y)$ such that
\begin{enumerate}[(i)]
\item for all $n_1, \ldots, n_k \in \NN$, $f(n_1, \ldots, n_k)$ is defined and equal to $m$\/ if and only if
$\vdash R_f(\numer{n_1}, \ldots, \numer{n_k},\numer{m})$
\item $\vdash \forall x: \synt{N}^k \forall y:\synt{N} \forall z : \synt{N}. R_f(x,y) \land R_f(x,z) \Rightarrow y=z$.
\end{enumerate}
\end{definition}

We now have the following characterization: 

\begin{theorem}[Lambek-Scott~\cite{LS86}] A partial numerical function is representable in HAH
(i.e., in the free topos) if and only if it is partial recursive.
\label{HAHparfns}
\end{theorem}

\subsection{C-monoids and Untyped Lambda Calculi}
 As mentioned earlier, Church's  untyped lambda calculus played a key role in 
 the original development of computability theory, as well as modern programming language
 theory.  It was Dana Scott in the late 1960s who pointed out that untyped lambda calculi may be 
 considered as typed lambda calculi with one non-trivial type (up to isomorphism). This arose from his development
 of {\em domain theory}, the mathematical modelling  of untyped lambda calculi and the semantics of
 programming languages \cite{AmadioR:domlc}. An algebraic framework for this development is given in \cite{LS86}, pp. 93-114, 
 which we now sketch.  For some historical references, the reader can see \cite{Lam80, DScott80}.

Recall, monoids are categories with one object.  A monoid has a terminal object precisely when it is trivial.
However, when we ignore the terminal object, we may formulate a notion of cartesian closure:
\begin{definition}[Lambek-Scott~\cite{LS86}]  A {\em C-monoid} is a monoid $\cM$ with constants $\pi_1,\pi_2, \varepsilon$, 
unary operation $(-)^*$, and binary operation $\langle - , - \rangle$ satisfying the equations of a CCC 
without a terminal object: i.e.  products, surjective pairing, $\beta$, $\eta$. Explicitly:

\begin{center}
$\pi_1\langle a, b\rangle = a \qquad \pi_2\langle a, b\rangle = b \qquad \langle \pi_1 c, \pi_2 c \rangle = c$ \\
$\varepsilon\langle h^*\pi_1,\pi_2\rangle = h \qquad (\varepsilon \langle k\pi_1,\pi_2\rangle)^* = k$
\end{center}
\end{definition}

The following results illustrate how C-monoids relate to untyped lambda calculi and CCCs.  They are 
an untyped variation of Theorem \ref{ccclam}.
\begin{theorem}[Lambek-Scott~\cite{LS86}]\
\begin{enumerate}[(i)]
\item There is a bijective correspondence between C-monoids and untyped lambda calculi with products and surjective pairing.
\footnote{Such untyped lambda calculi extended with surjective pairing do not enjoy good rewriting properties.  By a famous result of Klop
\cite{BarendregtH:lamcss}, Ch.15, \S 3, the Church-Rosser theorem fails for them.  Thus, the consistency of such systems would involve
constructing a non-trivial C-monoid (cf. \cite{LS86}, pp.107-114.) or more general models \cite{BarendregtH:lamcss}.  }
This correspondence extends to an isomorphism between the category of C-monoids and the category of such
untyped lambda calculi (cf. \cite{LS86}, p.106).
\item C-monoids correspond to CCC's generated by a non-trivial reflexive object $U$, i.e., an object  
$U \not \cong \synt{1}$ satisfying  $U^U\cong U \cong U\times U$.  In this case, $End(U)$ will be such a C-monoid
(cf. \cite{LS86}, p.99).
Without the \, $\eta$-rule, we would have  \, $U\times U \triangleleft U, \, U^U\triangleleft U$ (cf. also \cite{AmadioR:domlc},\cite{BarendregtH:lamcss} ).
\item With respect to appropriate numeral systems (e.g. Church or Barendregt numerals 
(see~\cite{BarendregtH:lamcss}, Sections 6.3, 6.4), the computable functions in the free C-monoid are precisely
the partial recursive functions (cf. \cite{LS86}, p.276.)
\end{enumerate}
\end{theorem}

We remark that part (ii) of the above theorem uses an observation of D. Scott (\cite{DScott80},\cite{LS86}), which says: if we form the 
idempotent splitting completion (Karoubi envelope) of a C-monoid, we obtain a CCC which is generated by a reflexive object $U$.  
There are precise senses in which all C-monoids are isomorphic to such CCCs (\cite{LS86}, p.99.)
Since Church's untyped lambda calculus was an early foundation of computability theory, 
it is no surprise that the computable functions in the free C-monoid are precisely the partial recursive ones.

\subsection{Plotkin's characterization of Kleene's $\mu$-recursion}

We recall Lambek's Lemma~\cite{Lam68a}, which is often used in denotational semantics.
 Given an endofunctor
$T: \cC\rightarrow \cC$ we define a {\em $T$-algebra} as  a map $TA\rightarrow A$.
Maps of $T$-algebras are commutative squares 
\[
\xymatrix{
TA \ar[r]^{Tf} \ar[d] & TB \ar[d] \\
A \ar[r]^f &B. }
\]
This gives a category of $T$-algebras; a $T$-algebra is called \emph{initial} when it is an initial object in 
this category.

\begin{lemma}[Lambek~\cite{Lam68a}]
If $TA\stt{\alpha} A$ is an initial $T$-algebra,  then $\alpha$ is an isomorphism.
\end{lemma}
For us, the following is the prime example:
\begin{example}{The NNO $\NN$ as an initial successor algebra in \set}
{Consider the endofunctor $T(-) = \synt{1} + (-)$ on $\set$\/ (often called the \emph{successor functor}), with the $T$-algebra
structure $(\synt{1} + \NN)\stt{\alpha} \NN$, where  $\alpha = [0,S]$, 
 for  $\synt{1}\stt{0} \NN$
and $\NN\stt{S}\NN$. The NNO property says that $\alpha$ is an initial $T$-algebra.
}
\end{example}

\noindent
In  $\set$,   Lambek's Lemma then gives the familiar fact that   
$\synt{1} + \NN \stt{\alpha} \NN$ is an isomorphism, for $\alpha = [0, S]$. 
 As we have seen above,  {\em initiality  of $\alpha$} gives us primitive recursion.
Now what about if we turn things around? Plotkin asked for the  {\em finality} of the co-algebra
$\alpha^{-1}: \NN \rightarrow \synt{1} + \NN$ -- not in  $\set$ but in  $\Par$.    
Interestingly, this turns out to give exactly Kleene $\mu$-recursion for partial functions.  

Let $\cC$ be a monoidal category with (right distributive) binary
sums and a weak left (or right) natural numbers object
$I \stt{0} \csynt{N} \stt{S} \csynt{N}$.  Following Plotkin, we extend Definition \ref{rep} of
{\em representable function} to include partial functions, as follows.
We shall say a partial function $f: \NN^k\rightharpoonup\NN$ is
{\em representable} by an arrow $F: \csynt{N}^k\rightarrow \csynt{N} \ \in \cC$ if
for all $n_1,\ldots,n_k\in \NN^k$,
\[
 f(n_1,\cdots,n_k)\cong m \tag{\dag} \ \Rightarrow \
 F\langle \numer{n_1}, \cdots , \numer{n_k}\rangle = \numer{m}
\]
where $\cong$ means Kleene equality.

\begin{theorem}[Plotkin~\cite{Plotkin13}]
Let $\cC$ be a monoidal category with (right distributive) binary
sums and a weak left (or right) natural numbers object
$I \stt{0}\csynt{N} \stt{S}\csynt{N}$ such that $[0, S]$ is an isomorphism and
$(\csynt{N}, [0, S]^{-1})$ is a weakly final natural numbers coalgebra.
Then all partial recursive functions are 
representable.
\end{theorem}

It is natural to ask if we can replace the ``$\Rightarrow$" in equation $(\dag)$ above 
 by the stronger condition ``$\Leftrightarrow$" (as in Definition \ref{parrep} (i))?
Plotkin calls this latter notion {\em strong representability}.  
The proof of Theorem \ref{HAHparfns} above (in \cite{LS86}, p.270) shows that 
for many arithmetical theories, representable partial functions are partial recursive. 
Plotkin takes the analog of this result (for strong representability) as an actual assumption to obtain a positive answer:
\begin{theorem}[Plotkin~\cite{Plotkin13}] Let $\cC$ be a monoidal category with (right distributive) binary
sums and a weak left (or right) natural numbers object
$I \stt{0} \csynt{N} \stt{S} \csynt{N}$ such that $[0, S]$ is an isomorphism and
$(\csynt{N}, [0, S]^{-1})$ is a weakly final natural numbers coalgebra. If
$0 \not = S0$ and if all strongly representable functions are partial recursive,
then all  partial recursive functions are strongly representable in $\cC$.
\end{theorem}

\section{Abstract Computability}  
 \label{sec-turing}
  In this section we address the question: what is a category of computable maps? This should be compared
 with ``synthetic" approaches to other areas of mathematics 
 such as synthetic differential geometry, synthetic domain theory, 
 homotopy type theory, and differential categories. A synthetic approach to computability 
 aims at describing the categorical structure
 common to all reasonable notions of computation; hence in such categories every morphism is by definition computable.
 Note the contrast with the work described in the previous section, where one starts with a category that, a priori, has
 no prescribed computability-theoretic content, and where one identifies some maps as representing 
 computable numerical functions.
 
 Most notions of computation are inherently partial, in the sense that they allow for the computable
maps to be partial maps. This fact, together with the importance of partial maps in other areas of mathematics,
has resulted in a long history of studying partial maps in categories, going back to the early days of topos theory.
This history largely overlaps with attempts to formulate aspects of computability theory in categorical terms, which
in turn are closely related to the study of categories of domains, as in \cite{AmadioR:domlc}.

\subsection{Categories of Partial Maps}
 
 We begin with a 
recent abstract treatment of categories of partial maps by Cockett and Lack~\cite{CockettJ:rescI}. There are at least 
two reasons for favouring this axiomatization: first, it is sufficiently general, in that it subsumes all
the previous treatments. Second, it is algebraic, in the sense that it identifies categories of partial maps as ordinary 
categories equipped with additional algebraic structure. This allows for the application of powerful techniques from
categorical algebra. For a much more detailed presentation and comparison with other approaches, see
loc. cit. and follow-ups.

\begin{definition}[Restriction Category] {A {\em restriction category} is a category $\cC$ together with
an assignment  $\overline{( \ )}: Hom_{\cC}(A,B) \longrightarrow Hom_{\cC}(A,A)$
mapping $f\longmapsto \overline{f}$ satisfying:
\vspace{1ex}

$\begin{array}{ll}
\mathrm{R.1}\quad &  f\overline{f} = f \\
\mathrm{R.2}\quad &  \overline{f}\overline{g} = \overline{g}\overline{f}$ whenever dom$\left(f\right) =$ dom$\left(g\right) \\
\mathrm{R.3}\quad &  \overline{g\overline{f}} = \overline{g}\overline{f}$ whenever dom$\left(f\right) = $ dom$\left(g\right) \\
\mathrm{R.4}\quad &  \overline{g}f = f\overline{gf}$ whenever cod$\left(f\right) = $ dom$\left(g\right) \\
\end{array}$
}
\end{definition}
We have $\overline{\overline{f}}=\overline{f}$, as well as ${\overline{f}} \, {\overline{f}}=\overline{f}$.
Maps satisfying $f=\overline{f}$\/ are called \emph{restriction idempotents}. The collection of restriction idempotents on
$A$\/ is denoted $\mathcal{O}(A)$; the composition operation makes $\mathcal{O}(A)$\/ into a meet-semilattice;
for each $f:A \to B$, there is an induced meet-semilattice homomorphism $f^*:\mathcal{O}(B) \to \mathcal{O}(A)$\/ 
sending $e \in \mathcal{O}(B)$\/ to $\overline{ef} \in \mathcal{O}(A)$.
A map 
$f: A\rightarrow B$ is {\em total} if $\overline{f} = id_A$. 
We obtain a wide subcategory $\mbox{\bf Tot}(\cC)\hookrightarrow \cC$.

\begin{example}{Examples of Restriction Categories}
\begin{enumerate}
\item $\Par$ is a restriction category when we define 
\[
\overline{f}(x) = \left\{
\begin{array}{ll}
x & \mbox{if $x\in Dom(f)$ }\\
\uparrow & \mbox{else}
\end{array}\right.
\]
\item The restriction structure on $\Par$\/ is inherited by various subcategories, most notably
the subcategory on the partial computable functions. 
This uses the fact that if $f$\/ is computable, then so is
$\overline{f}$.
\item Every category can be viewed as a restriction category by declaring $\overline{f}$\/ to be the identity for all $f$.
\end{enumerate}
\end{example} 

 \begin{definition}[Local Partial Order] For $f,g:A \to B$\/ in a restriction category,
 define
 \[ f \leq g \quad \Leftrightarrow \quad f = g\overline{f}.
 \] 
  \end{definition}   
\noindent
For example,  in $\Par$, we have:  $f\leq g$\/ precisely when $Graph(f) \subseteq Graph(g)$, i.e., when $g$\/ extends $f$.

Many notions for plain categories can be modified to make sense in the partial world. For example:
  
  \begin{definition}[Cartesian Structure]
A \emph{restriction terminal object}
is an object $\synt{1}$\/ together with, for each object $A$, a unique total map $!_A:A \to \synt{1}$\/ with 
$!_Bf \leq !_A$\/ for all $f:A \to B$. 

A \emph{restriction product} of $A,B$\/ is an object $A \times B$\/ 
with total projections $\pi_A,\pi_B$\/ such that for $f:C \to A, g:C \to B$\/ there is a unique 
$\langle f,g \rangle$\/ with $\pi_A\langle f,g \rangle \leq f$, $\pi_B\langle f,g \rangle \leq g$\/ and
$\overline{\langle f,g\rangle} = \overline{f}\,\overline{g}$. 
\[
\xymatrix{A \ar[r]^-{!_A} \ar[d]_f & \synt{1} \\
B \ar[ur]_{!_B}^\leq}
\qquad
\xymatrix{
& C \ar[dl]_f^\geq \ar@{->}[d]|{\langle f,g \rangle} \ar[dr]^g_\leq \\
A & A \times B \ar[r]_{\pi_B} \ar[l]^{\pi_A}  & B \\
}
\]
A {\em Cartesian Restriction Category} is a restriction category $\cC$ which 
has a restriction terminal object and binary restriction products. 
 \end{definition}

 \subsection{Turing Categories}\label{subsec:turing}
Turing categories, introduced in~\cite{CockettJ:T1} are restriction categories that essentially encode 
simultaneously the ideas underlying Kleene's $S^m_n$ and 
Enumeration theorems. They are also closely related to cartesian closed categories generated by
models of untyped lambda calculus, in that they weaken the cartesian closure, while generalizing to the partial world.
 
 \begin{definition}[Turing Category]
 A {\em Turing category} is a cartesian restriction category $\cC$
 with an object $A$\/ (called a \emph{Turing Object}),
  and a family of ``universal application morphisms" $\{A \times X \stt{\tau_{X,Y}}Y \ | \ X, Y \in \cC\}$ with 
  weak Currying:  for every $ Z\times X\stt{f}Y$ there exists a total map  $Z\stt{h} A$ factoring through $\tau_{X,Y}$:
\[
\xymatrix{
A \times X \ar[r]^-{\tau_{X,Y}} & Y \\
Z \times X \ar@{-->}[u]_{h \times 1} \ar[ur]_f \\
}
\]
 \end{definition}
Note that this expresses the idea that $A$\/ acts as a weak exponential $Y^X$, for any pair of objects $X,Y$.
One particular consequence is that every object is a retract of $A$.  In particular,
all finite restriction products $A^n$ are retracts of $A$. 

An elementary but useful fact is the fact that the class of Turing categories is closed under idempotent splitting:
if $\cC$\/ is a Turing category, then so is $\mathcal{K}_E(\cC)$\/ where $E$\/ is the class of restriction idempotents.

 A Turing category can equally well be described by ``universal self-application" 
 $\tau_{AA}$, denoted  $A\times A \stt{\bullet} A$.
 \begin{theorem}[Cockett-Hofstra~\cite{CockettJ:T1}] A Turing Category is a cartesian restriction category with
 an object $A$ such that (i) every object is a retract of $A$  and (ii) there is a universal
 self-application map $A\times A \stt{\bullet} A$.
 \end{theorem}
 Here are some of the motivating examples of Turing categories:
 
\begin{example}{Examples of Turing Categories}
\begin{enumerate}
\item Let $\{\phi_m\}_{m\in\NN}$ be a standard enumeration of unary partial recursive functions 
(see~\cite{CutlandN:comp}).
\emph{Kleene's First Model} $\comp(\NN)$ is the category whose objects are powers 
$\NN^k$ and whose maps $\NN^k\rightarrow \NN^m$ are $m$-tuples of partial computable 
functions of $k$ variables. $\NN$ is a Turing object,  there are retractions $\NN^k \triangleleft \NN$ and 
$m\bullet n := \phi_m(n)$ gives a universal application, by Kleene's theorems.
The restriction idempotents in this case are precisely the r.e. sets. Hence the restriction idempotent splitting of this 
category has the r.e. sets as objects, and partial computable functions as maps.

This example can be generalized to give categories $\comp(\NN^A)$, where $A$\/ is an oracle.
\item Consider a C-monoid, or more generally
 a reflexive object $U$ in a ccc, where $1\triangleleft U$, \ $U\times U \triangleleft U,$ \
 $U^U \triangleleft U$.  If  $ (m,r): U^U\triangleleft U$ is a retraction pair, then 
 $
 \bullet_U := U \times U \stt{r \times id_U} U^U \times U\stt{ev} U
 $
 determines a {\em total} Turing structure with Turing object U.
\item  Term models
 of  {\em Partial Combinatory Logic} (PCL) yield Turing categories. 
  PCL is an (partial) algebraic theory with constant symbols $\sss, \kkk$\/ and one binary application symbol 
  $\bullet$\/ (we write $xy$\/ instead of $x \bullet y$, and associate to the left). Terms are formed 
  in the usual way, together with a clause for forming restricted terms:
  \[ \mbox{\bf Terms} \qquad  t,t'::= \mathrm{VAR} \mid \sss \mid \kkk \mid tt'  \mid t|_{t'} \]
 where  $t |_{t'}$  is to be interpreted as ``$t$ restricted to $dom(t')$''. 
 (The categorical interpretation of such a restricted term
 is $\scbrl t \scbrr \circ \overline{\scbrl t' \scbrr}$.)
 The following equations are imposed: $\kkk xy=x_{|y}$, $\sss xyz=xz(yz)$, and $\sss xy \downarrow$. 
 See~\cite{CockettJ:T3} for details. 

The case of the \emph{closed term model}
is particularly significant because it corresponds to the initial Turing category.  Note that a total point $t:\synt{1} \to A$\/ 
of the Turing object corresponds to a provably total closed term $t$\/ of $PCL$. The global sections functor is therefore not
faithful, since there exist many closed terms that are not provably total, for example ${\kkk}|_{P}$\/ where $P$\/ is the
paradox combinator.
 \label{PCL}
 \end{enumerate}
  \end{example}
  
From the axioms of a Turing category, one may derive some basic results from computability theory such as the
recursion theorems. The restriction idempotents (partial identities serving as the domains of maps) in a Turing category
play the role of recursively enumerable sets; pullback of restriction idempotents then corresponds to m-reducibility.
Note that the standard model $\comp(\NN)$\/ also has \emph{ranges},
in the sense that every morphism not only has a domain but also a range; such categories are studied in detail
in~\cite{CockettJ:T4, CockettJ:T5}; see also~\cite{VinogradovaP:msc}.

Since the axioms of a Turing category 
do not preclude total models, one cannot expect results such as the undecidability
of the halting problem or Rice's theorem to follow in general. 
A detailed discussion of the development of basic computability
theory in the setting of Turing categories and how this depends on additional structure can be found in~\cite{Cockett10}.
 
\subsection{Computable maps and PCAs}
  
Turing categories are closely related to a class of structures called 
partial combinatory algebras (PCAs), as suggested by~Example \ref{PCL}. 
 
 Let $\cC$ be a cartesian restriction category. An \emph{applicative structure} in $\cC$\/ is a pair 
  $\cA = (A, \bullet)$, where $A\times A\stt{\bullet}A$\/ is a morphism called \emph{application}. There
  are no requirements on $\bullet$\/ (such as associativity). 
  Define $\bullet^n: A \times A^n \rightarrow A$ inductively, so $\bullet^{n+1} := 
  A \times A \times A^n\stt{\bullet \times id}A\times A^n\stt{\bullet^n}A$.

  \begin{definition}[Computable maps]   A map $A^n\stt{f}A$ is {\em $\cA$-computable}  when it is 
  ``named" by a total point of $A$, i.e. there is a total point $p: \synt{1}\rightarrow A$ such that   
  (identifying $A^n$ with $\synt{1} \times A^n$):
\[
\xymatrix{
A\times A^n \ar[r]^-{\bullet^n} & A   \\
     A^n \ar@{-->}[u]^{p \times id_A} \ar[ur]_f   \\
}
\]
  (Intuitively, $f(\vec{x}) = p\bullet \vec{x}$.)
Moreover, we require that $f$ is {\em total} on its first $n-1$ arguments.
More generally, we say a map $f:A^n\to A^m$ is {\em $\cA$-computable} if all its components are.
    \end{definition}
 
Since there are no axioms on an applicative object, the collection of $\cA$-computable maps cannot be expected to
have any good closure properties. In particular, it cannot be expected to form a subcategory of $\cC$. When it does,
the object $\cA$\/ is called \emph{combinatory complete}. 
This characterization is the categorical formulation of combinatory completeness (see
also~\cite{LongoG:cattcf}). 

Classically, an applicative structure is called combinatory complete
when every ``polynomial" built from variables, elements of $A$\/ and application, 
is represented by an element of $A$, see~\cite{BethkeI:notpca, OostenJ:reaics}. 
When $t$\/ is a polynomial and $x$\/ is a variable, we write
$\lambda^*x.t$\/ for the element representing $t$. That is: $(\lambda^*x.t)a=t[a/x]$\/ for all $a \in A$.
Equivalently, an applicative structure is a PCA exactly when it 
is a model of the theory PCL (see Example~\ref{PCL} above). 

\begin{definition} A combinatory complete applicative structure $\cA$ is called a \emph{partial combinatory algebra}
(PCA). For $\cA$\/ a PCA, denote by $\comp(\cA)$\/ the restriction category whose objects are the finite powers of $A$\/ 
and whose morphisms are the $\cA$-computable maps.
\end{definition}

At first sight, it may not be evident that combinatory completeness has many interesting consequences. It ensures, however,
that PCAs are powerful enough to represent various useful programming constructs. 

\begin{lemma}\label{lem:PCA}
In any PCA, we can define the following:
\begin{enumerate}[(i)]
\item Booleans, pairs, numerals (using e.g. $\numer{n}=\lambda^*a\lambda^*f.f^na$). 
What is more, any partial computable function $f:\NN \to \NN$\/ can then be represented in $\cA$\/ in the sense that
there is an element $a_f \in \cA$\/ for which $f(n)=m$\/ implies $a_f  \numer{n} = \numer{f(n)}$. 
\item Fixed points: there is an element $\synt{y} \in A$\/ for which $\synt{y}a = a(\synt{y}a)$\/ for all $a \in A$.
\item Recursors: there is an element $r \in A$\/ for which $rab\numer{0}=a, rab(\numer{n+1})=b(rab\numer{n})\numer{n}$\/ for
all $a,b \in A$\/ and $n \in \NN$.
\end{enumerate}
\end{lemma}

Standard examples of PCAs include Kleene's first model (natural numbers with partial recursive application, typically
denoted $\mathcal{K}_1$),
term models of PCL, and models of untyped lambda calculus. The following example is of importance in higher-order
computability, and will be used in the next section. We presuppose a surjective coding of finite sequences 
$\langle - \rangle:\NN^* \to \NN$. For $g:\NN \to \NN$, let $\hat{g}(k)=\langle g(0), \ldots, g(k-1)\rangle$. Finally,
let $*$\/ denote concatenation of sequences; for a sequence $L$\/ and $n \in \NN$\/ we write $n * L$\/ instead
of $\langle n \rangle * L$.

\begin{example}{Kleene's second model}\label{ex:K2}
Consider $f,g \in \NN^\NN$. Define a (possibly partial) function $f \star g:\NN \to \NN$\/ by
\[
(f \star g)(n) = \left\{ \begin{array}{ll} 
f(\langle n*\hat{g}(k) \rangle)-1 & \textrm{ where } k=\mu_r.f(\langle n*\hat{g}(r) \rangle)>0 \\
\textrm{undefined} & \textrm{ if no such $k$\/ exists}. \\
\end{array}
\right.
\]
Then define a partial application $\bullet:\NN^\NN \times \NN^\NN \to \NN^\NN$\/ by
\[ f \bullet g =\left\{ \begin{array}{ll}  
f \star g & \textrm{ if } f \star g \textrm{ is total } \\
\textrm{undefined} & \textrm{ otherwise}.
\end{array} 
\right.
\]
This model is typically denoted $\mathcal{K}_2$, and captures a notion of ``computable operations acting on
continuous data". If we restrict $\NN^\NN$\/ to the set of total computable functions, we get a 
sub-PCA $\mathcal{K}_2^{\mathrm{eff}}$\/ 
of ``computable operations acting on computable data''. 
For details see~\cite{LongNorm15}.
\end{example}

What is the correct notion of morphism of PCAs? Regarding $\cA, \cB$\/ as computational devices, 
a morphism $\varphi:\cA \to \cB$\/ should at least express that $\cA$\/ can be interpreted, or simulated, within $\cB$.
The following definition is due to Longley~\cite{LongleyJ:reatls}. We state it in set-theoretic terms here, but it
can easily be rendered diagrammatically in a cartesian restriction category:

\begin{definition}[Simulation]
A \emph{simulation} from $\cA$\/ to $\cB$\/ is a function $\varphi:A \to B$\/ for which there exists $b \in B$\/ such that
\[ x\bullet y\!\downarrow \quad \Longrightarrow \quad b \bullet \varphi(x) \bullet \varphi(y) = \varphi(x\bullet y).\]
\end{definition}
Simulations compose, and in fact form a 2-category. 
We point out that in~\cite{LongleyJ:reatls} a relational version of this definition is given; however, as demonstrated 
in~\cite{HofstraP:ordpca}, it is possible to view relational simulations as Kleisli morphisms over a base category
of functional simulations.

\begin{example}{Numerals as Simulation}
\noindent
Every PCA admits a choice of numerals; such a choice amounts to a simulation $\mathcal{K}_1 \to \cA$. 
All non-constant simulations $\mathcal{K}_1 \to \cA$\/ are in fact isomorphic to each other.
\end{example}

 Returning to the connections between PCAs and Turing categories, we 
 note that by construction $\comp(\cA)$\/ is a cartesian restriction category. 
  The following shows that PCAs are a fundamental notion for generating
 Turing categories: every PCA gives rise to a Turing category, and every Turing category is generated by
 the PCA structure on the Turing object.
 \begin{theorem}[Cockett-Hofstra~\cite{CockettJ:T1}]~~~~~
 \begin{enumerate}[(i)]
\item If $\cA = (A, \bullet)$ is a PCA, then $\comp(\cA)$ is a Turing category, with Turing object $A$.
\item If $(\cC,A)$ is a Turing category with Turing object $A$, then $(A,\bullet)$ is a PCA and 
$\cC \cong \mathcal{K}_{E}(\comp(\cA))$, for some class of idempotents $E$.
\end{enumerate}

 \end{theorem}
   
Thus  ``{\em Categories of the form $\comp(\cA)$ serve as a minimal environment (for) 
PCA's and ... computable maps ...;  other Turing categories are supposed to be viewed as (non-essential) 
inflations of such minimal categories}" (\cite{CockettJ:T1}).

Earlier we contrasted the approach of identifying representable numerical functions in free categories with NNO with
the synthetic approach of Turing categories. However, there is a slightly different perspective on Turing categories, that perhaps
brings the two approaches closer together.

Instead of considering Turing categories in isolation, i.e., synthetically, one can consider Turing categories
structured over a base category. For example, the Turing category
$\comp(\NN)$\/ can be considered as a non-full subcategory of $\Par$. This point of view is particularly relevant when one
wishes to consider non-computable functions or study, e.g., non-r.e.degrees.
More generally, we think of a Turing category $\cC$\/ 
with a cartesian restriction functor $F:\cC \to \cB$\/ into a base category $\cB$\/ as specifying an object $FA$\/ of $\cB$
together with a notion of computation on $FA$. The object $FA$\/ is necessarily a PCA, but $\cC$\/ is not
always $\comp(FA)$; the reason is that $FA$\/ may have total elements $t:\synt{1} \to FA$\/ that are not in the image of $F$.
Hence $\comp(FA)$\/ may contain morphisms that are not represented in $\cC$. This forces the consideration of
\emph{relative PCAs}, and the full characterization of Turing categories over a fixed base in terms of such 
relative PCAs can be found in~\cite{CockettJ:T2}.\footnote{This characterization involves a notion of \emph{simulation}
between Turing categories (over a fixed base), generalizing the foundational work by Longley~\cite{LongleyJ:reatls} on
simulations between PCAs (called \emph{applicative morphisms} in loc. cit.).}

Note that there is an analogy between the two perspectives on Turing categories and those on toposes: 
one may consider toposes relative to a fixed base topos $\mathcal{S}$\/ (as is common in the study of Grothendieck toposes,
where $\mathcal{S}$\/ plays the role of the universe of sets), or one may study elementary toposes such as the free topos
without regarding them as being constructed over a base.

 
\section{Realizability}\label{sec:realizability}
We now briefly turn our attention to a strand of research that also heavily involves the study of categorical structures
associated to models of computation, but that is different from the earlier themes in that it primarily considers such structures
as models of various logical systems. 

\subsection{Kleene Realizability}
Realizability, originally devised by Kleene in the seminal paper~\cite{KleeneS:intint}\footnote{We omit a discussion of
the history of the subject, of which some of the main threads are detailed in~\cite{OostenJ:reahe}.}, 
is to be thought of as
a semantics for constructive mathematical systems\footnote{Recent work by Krivine and others has shown that it
is also possible to define realizability interpretations of classical systems.}. In Kleene's original work, the system at hand
was Heyting Arithmetic (HA), and the interpretation was defined in terms of partial computable functions. The central notion
is written $\real{n}{A}$, where $n \in \mathbb{N}$\/ and $A$\/ a formula in the language of arithmetic, and is pronounced
``$n$\/ realizes $A$", or ``$n$\/ is a realizer for $A$". The intuition is that $n$\/ codes information about why $A$\/ is
true. The definition is by induction on the structure of $A$\/ (and uses an enumeration $\phi_0, \phi_1, \ldots$\/) of
unary partial computable functions:

\begin{definition}[Kleene Realizability]
Define $\real{n}{A}$\/ (for sentences $A$) by

\begin{tabular}{lcl}
$\real{n}{t=s}$\/ & iff & $n=0$\/ and $t=s$\/ is true \\
$\real{n}{A \land B}$\/ & iff & $n=\langle a,b \rangle$\/ where 
$\real{a}{A}$\/ and $\real{b}{B}$ \\
$\real{n}{A \lor B}$\/ & iff &  $n=\langle a,b \rangle$\/ where either
$a=0$\/ and $\real{b}{A}$\/ or $a=1$\/ and $\real{b}{B}$ \\
$\real{n}{A \to B}$\/ & iff & for all $m \in \NN$, 
if $\real{m}{A}$\/ then $\phi_n(m) \!\downarrow$\/ and  $\real{\phi_n(m)}{B}$ \\
$\real{n}{\exists x.A}$\/ & iff &$n=\langle a,b \rangle$\/ where 
$\real{b}{A[a/x]}$ \\
$\real{n}{\forall x.A}$\/ & iff & for all $m \in \NN$, 
$\phi_n(m) \!\downarrow$\/ and $\real{\phi_n(m)}{A[m/x]}$ \\
\end{tabular}
\end{definition}

The Soundness theorem now states: $HA \vdash A \;  \Longrightarrow \; \exists n \in \NN. \real{n}{A}$. The converse,
however, is false: there are realizable statements that are underivable.
Most notably, \emph{Extended Church's Thesis} ($ECT_0$)
is the scheme:
\[
\forall x(A(x) \to \exists y.B(x,y)) \to \exists e\forall x(A(x) \to B(x,e \bullet x))
\]
Here, $A$\/ is assumed to be an almost negative formula, and 
$e \bullet x$\/ denotes the application of the $e$-th computable function to $x$, suitably represented in $HA$.
One can show that all instances of $ECT_0$\/ are realizable but not provable in $HA$. Moreover, $ECT_0$\/ 
\emph{axiomatizes} Kleene realizability, in the sense that the realizable statements of $HA$\/ are precisely those
that are derivable in $HA+ECT_0$.

Over the years, many variations on Kleene's original definition have been studied, with the purpose of establishing,
among other things, consistency results and proof-theoretic properties of various formal systems. For example,
\emph{q-realizability} incorporates derivability into the definition of realizability, and can be used to establish
the existence and disjunction properties of HA. 
     
\subsection{Realizability Toposes}
How does realizability manifest itself categorically? Historically, the topos-theoretic treatment of
Boolean-valued and Heyting-valued models (\cite{FourmanM:shel, HiggsD:catabv}) inspired the idea of considering 
sets of realizers as truth values. This idea led Hyland to the discovery of the \emph{Effective Topos}~\cite{HylandJ:efft}, 
an elementary (non-Grothendieck) topos with the property that the first-order arithmetical statements about the NNO 
are precisely the Kleene-realizable statements. Thus, among other things, the internal language of $\eff$\/ is a 
natural extension of Kleene realizability to higher-order logic. 

Various notions from computability theory find a natural home in $\eff$. 
For example, the Turing degrees manifest themselves in the form of subtoposes of $\eff$:

\begin{theorem}[\cite{HylandJ:efft, PhoaW:relet}]
The lattice of Turing degrees faithfully embeds into the lattice of subtoposes of $\eff$.
\end{theorem}

(Here, the notion of subtopos is taken in the geometric sense: it is a full subcategory closed under finite limits, whose inclusion
has a finite-limit preserving left adjoint.) Not every subtopos arises from a Turing degree however; 
see~\cite{LeeS:basse} for more information. 

There are several ways to present the Effective Topos and its variants. Perhaps the simplest is via \emph{exact completions}
(see~\cite{CarboniA:somfcr, CarboniA:regec}, as well as~\cite{MenniM:chalex}). A category $\cC$\/ is called \emph{exact} if it has finite limits,
pullback-stable quotients of equivalence relations, 
and if every coequalizer is the coequalizer of its kernel pair. Every topos is exact.
Now to each category with finite limits $\cC$\/ one may associate an exact category $\cC_{ex/lex}$\/ by freely adding quotients
of equivalence relations, and the Effective Topos is of this form. 
The finite limit category in question is called $\pass$, the category of 
\emph{partitioned assemblies}.

\begin{definition}[Partitioned Assemblies]
The category $\pass$\/ has objects pairs $(X,\alpha)$\/ with $X$\/ a set and $\alpha:X \to \NN$\/ a function; a morphism
$(X,\alpha) \to (Y,\beta)$\/ is a function $f:X \to Y$\/ which is \emph{tracked}, in the sense that there exists a code $e \in \NN$
such that 
\[ \forall x \in X. e \bullet \alpha(x)\! \downarrow \; \land \ e \bullet \alpha(x) = \beta(f(x)).\]
\end{definition}

\begin{theorem}[Carboni et al.~\cite{CarboniA:somfcr, CarboniA:catarp}]
The Effective Topos is the exact completion of the category of partitioned assemblies: $\eff \simeq \pass_{ex/lex}$.
\end{theorem}

The above construction of $\eff$\/ can be refined by considering an intermediate category:
\begin{definition}[Assemblies]
The category $\ass$\/ has objects pairs $(X,\alpha)$\/ with $X$\/ a set and $\alpha:X \to \mathcal{P}_+\NN$\/ a function
(where $\mathcal{P}_+$\/ denotes the non-empty powerset); 
a morphism $(X,\alpha) \to (Y,\beta)$\/ is a function $f:X \to Y$\/ which is \emph{tracked}, 
in the sense that there exists a code $e \in \NN$ such that 
\[ \forall x \in X \forall a \in \alpha(x). e \bullet a\! \downarrow \; \land \ e \bullet a \in \beta(f(x)).\]
\end{definition}

The category $\ass$\/ is \emph{regular}: it has finite limits and admits stable quotients of equivalence relations. 
Any finite limit category $\cC$\/ admits a \emph{free regular completion} $\cC_{reg}$, and any regular category $\cD$\/ 
admits a free exact completion $\cD_{ex/reg}$. 
With this notation, we now have the following relations between $\pass, \ass$, and $\eff$:

\begin{theorem}[Carboni et al.~\cite{CarboniA:somfcr, CarboniA:catarp}]
There are equivalences $\ass\simeq \pass_{reg}$\/ and $\eff \simeq \ass_{ex/reg}$.
\end{theorem}
The category of assemblies happens to be much more than regular: it is a quasitopos and has a NNO,
given by $(\NN,\{-\})$. As such, a lot of the computability-theoretic features of $\eff$\/ already manifest themselves in 
this subcategory. For example, in $\ass$\/ we may consider higher-type computability over $\NN$.

An alternative construction of $\eff$, more logical in nature, makes use of the concept of a \emph{tripos} 
(see~\cite{HylandJ:trit}; tripos is an acronym for ``topos-representing indexed preordered set".) One considers
the $\set$-indexed preorder $\set(-,\mathcal{P}\NN)$; for a set $X$, we preorder $\set(X,\mathcal{P}\NN)$\/ by:
\[ \alpha \vdash_X \beta \Longleftrightarrow \exists e \in \NN \forall x \in X \forall a \in \alpha(x). e \bullet a \!
\downarrow \; \land \ e \bullet a \in \beta f(x).\]
There is now a general construction turning a tripos into a topos, and $\eff$\/ arises as the result of applying this
construction to $\set(-,\mathcal{P}\NN)$. This construction highlights the original idea of regarding sets of realizers
as truth-values, in analogy with $H$-valued sets for $H$\/ a complete Heyting algebra.

\subsection{PCAs and Toposes}
The construction of the Effective Topos generalizes in various ways. We focus on the following fact\footnote{It was 
already known well before the discovery of the effective topos
 that combinatory algebras carried sufficient structure to define notions of realizability,
 see e.g.~\cite{FefermanS:lanaem}.}: for each
PCA $\cA=(A,\bullet)$, there is an associated realizability topos $\RT{\cA}$. In fact, we may associate to $\cA$\/ 
a category of partitioned assemblies $\pass(\cA)$\/ (where the objects are sets $X$\/ equipped with
a function $\alpha:X \to A$), and let $\RT{\cA}=\pass(\cA)_{ex/lex}$. Alternatively
we build the tripos $\set(-,\mathcal{P}A)$. The functoriality of $\cA \mapsto \RT{\cA}$, including the correct notion
of ``Morita equivalence" for PCAs was worked out in~\cite{LongleyJ:reatls}; the complete 
characterization of (geometric) morphisms between toposes of the form $\RT{\cA}$\/ in terms of
morphisms of (ordered) PCAs appears in~\cite{HofstraP:ordpca}.

An important construction, both for the analysis of realizability toposes and for applications of realizability, is
that of the category of PERs over a PCA. A PER (partial equivalence relation) on a set $A$\/ is simply a symmetric
and transitive relation; equivalently, it is an equivalence relation on a subset of $A$\/ (then called the \emph{domain}
of the PER). When $R$\/ is a PER on $A$, we write $A/R=\{[a] \mid (a,a) \in R\}$\/ for the set of equivalence classes.
In case of a PCA, this leads to the following:

\begin{definition}[Category of PERs]
Let $\cA=(A,\bullet)$ be a PCA. The category $\per{\cA}$\/ has as objects PERs $(A,R)$\/ on $A$. A morphism 
$(A,R) \to (A,S)$\/ is a function $f:A/R \to A/S$\/ that is tracked in the sense that there exists $a \in A$\/ such that
\[ \forall x \in A. (x,x) \in R \to a \bullet x\!\downarrow \land \; f[x]=[a \bullet x].\]
\end{definition}

The category $\per{\cA}$\/ can be seen as a full subcategory of $\ass(\cA)$\/ on those objects $(X, \alpha)$\/ for which
$\alpha(x) \cap \alpha(y) \neq \emptyset$\/ implies $x=y$. It is (locally) cartesian closed, and has a NNO. We will 
return to this structure in the section on higher type computability below.

Since PCAs give rise both to Turing categories and to realizability toposes, it is natural to wonder how the latter two are
related. We mention here one result that builds on earlier insights into how realizability toposes
can be regarded as colimit completions~\cite{RobinsonE:colcet, RobinsonE:abslr}. In~\cite{CockettJ:T2} 
 a universal property of partitioned assemblies is exhibited: it is the free fibred preorder on a functor,
in a suitable restriction-category theoretic sense. In case of a PCA $\cA$\/ with associated Turing category 
$\comp(\cA)$, applying this construction gives a fibration, and taking total maps recovers $\pass(\cA)$.
Moreover, this construction has the property that it turns simulations
between Turing categories into actual functors on the level of partitioned assemblies. 

To conclude our discussion of realizability toposes  we mention the abstract characterization of toposes of the form $\RT{\cA}$\/ 
due to Frey~\cite{FreyJ:chapar}. In order to state this result, we need to define a few concepts. First, suppose that
$\Gamma \dashv \nabla$\/ is a pair of adjoint functors with $\Gamma \circ \nabla \cong 1$. Then a map 
$f$\/ is called \emph{closed} (w.r.t. this adjunction) if the square
\[
\xymatrix{
A \ar[r]^f \ar[d] & B \ar[d] \\
\nabla\Gamma A \ar[r]_{\nabla \Gamma f} & \nabla\Gamma B \\
}
\]
in which the vertical maps are the unit morphisms is a pullback\footnote{The terminology \emph{closed} derives from 
the fact that for realizability toposes $\RT{\cA}$, closed subobjects for the double negation topology are characterized
by this condition.} . Moreover, an object $A$\/ is called \emph{separated}
when the unit $A \to \nabla\Gamma A$\/ is monic\footnote{This terminology also derives from the fact that in $\RT{\cA}$\/ 
this characterizes the separated objects for the double negation topology.}. 
Finally, an object is called \emph{discrete} when it is orthogonal to all
closed regular epimorphisms.

\begin{theorem}[Frey~\cite{FreyJ:chapar}]
A locally small category $\cE$\/ is equivalent to $\RT{\cA}$\/ for a PCA $\cA$\/ if and only if the following conditions hold:
\begin{itemize}
\item $\cE$\/ is exact and locally cartesian closed;
\item $\cE$\/ has enough projectives and the full subcategory $\mathrm{Proj}(\cE)$\/ on the projective objects is closed
under finite limits;
\item The global sections function $\Gamma:\cE \to \set$\/ has a right adjoint $\nabla$\/ which factors through 
$\mathrm{Proj}(\cE)$;
\item There exists a separated, projective object $D$\/ such that for any projective object $P$\/ there exists a closed
map $P \to D$.
\end{itemize}
\end{theorem}

This theorem should be regarded as the analogue of the well-known Giraud theorem characterizing Grothendieck toposes 
among exact categories in terms of their relation to $\set$. Note that the first conditions express that $\cE$\/ is of the form
$\cC_{ex/lex}$, and that the last two conditions therefore characterize categories of the form $\pass(\cA)$.

We end this section by a brief mention of another approach to partial recursive functions and PER, introduced by Lambek
\cite{Lambek97} and studied further in \cite{LambekJ:exampr}.  In this view, one considers the  category of relations generated by the 
monoid of primitive recursive functions (qua relations). Taking this viewpoint, a partial recursive function is simply a single-valued recursively enumerable (r.e.) relation, and the category PER is a kind of Karoubi envelope construction:  
the category whose objects are arbitrary pers on $\NN$ and whose maps are r.e. functional relations 
between them.  The full subcategory of PER given by r.e. pers and r.e. functional relations 
is particularly interesting in this regard, since it turns out to be exact. In 
\cite{LambekJ:exampr}, it is considered as a candidate for a kind of 
exact completion of the monoid of primitive recursive functions, although the precise nature of this completion is yet
to be determined.

\section{Other Directions}
This final section briefly introduces some facets of computation that have a somewhat different character than the
work discussed so far. First, we discuss traced monoidal categories and PCAs arising in ``Geometry of Interaction" situations.
Next, we turn to computability at higher type, giving a very brief introduction to some of the concepts and ideas in that area.
Finally, we mention some of the categorical approaches to complexity theory.

\subsection{Traced Categories}\label{subsec:trace}
In an influential paper, Joyal, Street, and Verity \cite{JSV96} introduced the notion of an abstract trace in monoidal categories.
Such traces arise in a wide range of areas, including knot theory, fixed point theory and theoretical computer science.  
We will be especially concerned with  applications arising in the algebra of feedback in networks and the associated 
fixed point theories. Traced monoidal categories also play a prominent
role in the categorical analysis of Girard's Geometry of Interaction (GoI) Program in Linear Logic, in which one  
analyzes the dynamics and flow of information in cut-elimination in networks of proofs~\cite{AHS02,HS06}.  
For simplicity, we consider the case of symmetric monoidal categories.

A {\em parametrized trace} on a symmetric monoidal category $\cC$ is an operation $tr:  \cC(X\otimes U, Y\otimes U)
 \rightarrow \cC(X,Y)$, satisfying a number of axioms discussed in detail in \cite{JSV96,AHS02}. 
 The theory has a particularly geometric
 flavour, and the papers, loc. cit., use a {\em string calculus} both for describing the axioms and for diagrammatic reasoning.   
 A particular evocative picture is to think of the trace as a form of ``feedback":
 
\vspace{1cm}
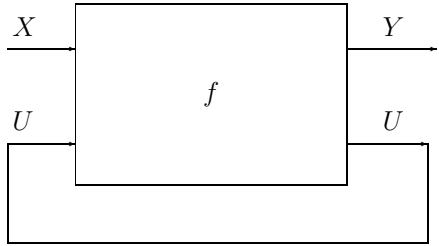
 \[
 \scalebox{.6}{
 \begin{picture}(100,30)(5,0)
\put(5,5){\framebox(60, 40){\huge  $f$}}
\put(-9, 38){\huge $ X$}
\put(-9, 17){\huge $ U$}
\put(73, 38){\huge $ Y$}
\put(73, 17){\huge $ U$}
\put(65, 35){\vector(1,0){20}}
\put(65, 14){\vector(1,0){18}}
 \put(-10, 35){\vector(1,0){15}}
\put(-10, 14){\vector(1,0){15}}
\put(-10, 14){\line(0,-1){22}}
\put(83, 14){\line(0,-1){22}}
\put(-10, -8){\line(1,0){93}}
\end{picture}
}
 \]

\vspace{.5cm}
\noindent         
Examples relevant to this paper include \rel \ and \Par, with $\otimes = \biguplus$, the disjoint union of sets.  
In the case of \Par, the trace of a map 
$f: X\uplus U\rightarrow Y\uplus U$ is given by the  following summation formula:
\[
Tr^{U}_{X,Y}(f) = f_{XY} + \sum_{n\in \NN}f_{UY}f^n_{UU}f_{XU}
\]
Here $f_{XY}$ denotes the partial map $X\rightharpoonup Y$ obtained from $f$ by naturally restricting 
the domain and codomain (using injections and partial projections), and similarly for the other components.  
The sum of a family of partial maps $\sum_{n\in \NN}h_n: X\rightharpoonup Y$ is defined iff the domains of the $h_n$ are disjoint, 
in which case $(\sum_{n\in \NN}h_n)(x) = h_k(x)$ if $h_{k}(x)\!\!\downarrow$, and is undefined otherwise. 
Such traces given by the above formula are called ``particle-style" (\cite{AHS02}) based on the following intuition:  
in the above picture imagine particles 
entering the box at $X$.  Either they exit immediately at $Y$\/ via $f_{XY}$\/
or they exit through $U$\/ and continue to cycle on $U$ some finite number $n$\/ times via $f_{UU}$\/ and 
then eventually exit at $Y$.

In~\cite{AHS02}, it is shown how a so-called \emph{GoI situation} gives rise to a 
\emph{linear combinatory algebra}. 
A GoI situation is a traced symmetric monoidal category equipped with a traced symmetric monoidal endofunctor,
and an object $U$\/ satisfying various domain equations. By applying the \emph{GoI construction}, one obtains
a compact closed category containing an object whose points form a linear combinatory algebra. By the latter,
one means an applicative structure $(A,\bullet)$\/ equipped with an endomap $!:A \to A$\/ and several combinators, allowing 
for the application $(a,b) \mapsto a \bullet !b$\/ to form a total combinatory algebra.

Lambek's register machines were described by a language of flowcharts and feedback.  
They can be naturally represented in a symmetric traced category 
with $\otimes$ = coproduct~\cite{KatSabWalt02}.  The original categorical studies of  iterative notions of flowchart 
computation in a programming language setting was by C. Elgot.   
In this case iteration is given by a kind of feedback loop in a category whose hom-sets have infinite sums 
(Elgot's ideas are detailed in~\cite{ManArb86}, and pursued more abstractly in traced $\Sigma$-monoid enriched tensor 
categories by Haghverdi~\cite{Hagh00}).  Finally, traced monoidal categories 
in which the monoidal tensor $\otimes$ is obtained from a cartesian or genuine tensor product are discussed 
in~\cite{AHS02}, as well as a more
general notion of partially traced categories, in~\cite{HS10}.

\subsection{Typed PCAs}
The notions of computation considered so far has been \emph{untyped}, in the sense that it is based on a 
single base type containing both the input/output values of computable maps and the (codes for) computable maps.
In the notion of PCA, this is reflected by the fact that the partial application $a \bullet b$\/ regards $a$\/ as 
a code for a partial map and $b$\/ as an input. In various situations however, we do wish to consider computation
over different types, for example because we wish to distinguish between the type of computable operations and the type
of its inputs and outputs. One of the key concepts in the study of such situations is that of a typed PCA.

\begin{definition}[Typed PCA]
Let $T$\/ be the collection of simple types generated by a single base type $N$.  
A \emph{Typed Partial Combinatory Algebra} (TPCA) over $T$\/ is a set-valued assignment
$\tau \mapsto A(\tau)$\/ for $\tau \in T$, together with for all $\sigma, \tau \in T$\/ a partial application function 
$\bullet_{\sigma,\tau}: A(\sigma \to \tau) \times A(\sigma) \to A(\tau)$. As for PCAs, 
we write application using infix notation, associating to the left; we also suppress the typing information.
 One requires the existence of 
combinators 
\[ {\kkk}_{\sigma,\tau} \in A(\sigma \to (\tau \to \sigma)) \; ; \qquad {\sss}_{\sigma, \tau, \rho} \in 
A((\sigma \to \tau \to \rho) \to ((\sigma \to \tau) \to (\sigma \to \rho)))\]
 (for all types $\sigma, \tau, \rho$)
satisfying
\begin{itemize}
\item ${\kkk} \bullet x \bullet y = x$
\item ${\sss} \bullet x \bullet y\!\downarrow$ 
\item ${\sss} \bullet x \bullet y \bullet z = (x \bullet z) \bullet (y \bullet z)$
\end{itemize}
\end{definition}

We remark that some authors also require the existence of \emph{fixed point combinators}, \emph{numerals}
and \emph{recursors} (see Lemma~\ref{lem:PCA} for what these are in the untyped setting). 
This essentially guarantees that a TPCA is a model of Plotkin's simply typed programming 
language PCF (see~\cite{PlotkinG:lcfcpl}).

\begin{example}{Examples of TPCAs}
\begin{enumerate}
\item Let $A(N)=\mathbb{N}$, and $A(\sigma \to \tau)=A(\tau)^{A(\sigma)}$. Then
we can let application be evaluation $\bullet:A(\tau)^{A(\sigma)} \times A(\sigma) \to A(\tau)$. This is called
the \emph{full} (total) TPCA over $\mathbb{N}$.
\item In the previous example we may instead let $A(\sigma \to \tau)=\Par(A(\sigma),A(\tau))$, the set of all 
partial functions. Then we get a TPCA where application is partial.
\item If $\cC$\/ is a CCC with NNO, we consider the subcategory on the simple types over the NNO. Taking
global sections gives a TPCA.
\item Any PCA $\cA=(A,\bullet)$\/ is a typed PCA where 
$A(\sigma)=A$, and $\bullet_{\sigma,\tau}=\bullet$\/ for all types $\sigma, \tau$.
\item Term models of typed lambda calculus form TPCAs in the expected way, as do term models of
programming languages based on typed lambda calculus, such as PCF.
\end{enumerate}
\end{example}

Just as for PCAs, there is a notion of \emph{simulation} between TPCAs. For example, to say that 
$\cA$\/ has numerals (and that computable functions are representable) is to say that there is a simulation 
of Kleene's first model into $\cA$. See~\cite{LongNorm15} for details.

\subsection{Computation at higher types}
Most of the developments described above concern first-order computability (possibly taking place in a higher-order
setting). We now briefly discuss computability at higher types. The relation between higher-order computability and
first-order computability is analogous to that between functional analysis and analysis. Thus in higher-order computability
one studies functionals $\NN^\NN \to \NN$, and so on. Immediately, one recognizes the many possibilities: one could
consider functionals acting on all total functions, or on all partial functions, or on all total computable functions, or on all 
partial computable functions, et cetera. We refer to the detailed survey paper~\cite{LongleyJ:notchtI} for a 
comprehensive historical overview. 

\begin{example}{Kleene's S1-S9}
One of the most fundamental notions of higher type computability was introduced in the landmark paper~\cite{KleeneS:recfqf}.
The collection of \emph{pure} types over $\NN$\/ is defined by: 
\[ \NN^0 = \synt{1}, \quad \NN^{(k+1)} = \NN^{\NN^{(k)}}. \]
Kleene's conditions S1-S9 define a class of partial maps of type
\[ \Phi: \NN^{(k_1)} \times \cdots \times \NN^{(k_r)} \to \NN. \]
More precisely, the definition specifies a relation $\{e\}(v_1, \ldots, v_r) \simeq x$, where $e \in \NN$\/ is an index,
the $v_i$\/ are elements
of the pure types $\NN^{(k_i)}$, and $x \in \NN$. Thus the resulting definition is an example of partial functionals operating
on total functions. 
\end{example}

Another classic example of a notion of computation at higher type, first introduced in~\cite{KreiselG:intamc}, is the following:

\begin{example}{Hereditarily Effective Operations}
Define simultaneously, for each simple type over the natural numbers, a set of natural numbers
and an equivalence relation on the set as follows:
\begin{itemize} 
\item $\HEO_0=\NN$, and $n \sim_0 m \Leftrightarrow n=m$.
\item $\HEO_{\sigma \to \tau}=\{e \in \NN \mid \phi_e$\/  induces a total function $\HEO_\sigma \to \HEO_\tau\}$,
and $e \sim_{\sigma \to \tau} e' \Leftrightarrow \forall n \in \HEO_\sigma. \phi_e(n) \sim_{\tau} \phi_{e'}(n)$.
\end{itemize}
\end{example}

One of the central contributions in~\cite{LongNorm15} is the development of a general framework (called 
\emph{computability model})
 for studying the wide variety of possible notions
of higher type computation. It also supports a general notion of simulation between models, and of equivalence. 
Typed PCAs form an important class of examples of computability models. We shall now sketch a result
by Longley characterizing the so-called extensional collapse of a large family of TPCAs. From now on, we assume
our TPCAs come equipped with a choice of numerals $\mathbb{N} \to A(N)$. 

\begin{definition}[Extensional Collapse of a TPCA]
Let $\cA$\/ be a TPCA. Define, at each simple type $\sigma$, a PER $\sim_\sigma$\/ on $A(\sigma)$\/ as follows:
\begin{itemize}
\item $a \sim_N b$\/ iff $a=b=\overline{n}$\/ for some $n \in \NN$
\item $a \sim_{\sigma \to \tau} b$\/ iff for all $x,y \in A(\sigma)$\/ with $x \sim_\sigma y$, $a \bullet x \sim_\tau b \bullet y$.
\end{itemize}
The sets $A(\sigma)/\!\sim_\sigma$\/ form a simple type structure over $\NN$, denoted $\EC{\cA}$.
\end{definition}

\begin{definition}
A typed PCA $\cA$\/ is 
\begin{enumerate}[(i)]
\item \emph{continuous} if there is a numeral-respecting simulation $\cA \to \mathcal{K}_2$;
\item \emph{full continuous} if it is continuous and all functions $\mathbb{N} \to \mathbb{N}$\/ are
represented in $\cA$; 
\item \emph{effective} if there is a numeral-respecting simulation $\cA \to \mathcal{K}_1$; 
\end{enumerate}
\end{definition}

The following general result (referred to as the \emph{Ubiquity Theorem}) now describes the extensional collapse of
these important classes of typed PCAs:

\begin{theorem}[Longley~\cite{LongleyJ:ubictt}]
Let $\cA$\/ be a typed PCA.
\begin{enumerate}[(i)]
\item If $\cA$\/ is full continuous, then $\EC{\cA}=\mathsf{C}$, the total continuous functionals
(which may be taken to be $\mathsf{C}=\EC{\mathcal{K}_2}$).
\item If $\cA$\/ is effective (and satisfies a few minor technical conditions), then $\EC{\cA}=\HEO$, the 
hereditarily effective operations.
\end{enumerate}
\end{theorem}

There is a third part to the theorem, which characterizes the collapse of a class of \emph{relative} TPCAs. By the latter,
we mean a TPCA $\cA$\/ together with a sub-TPCA $\cA^\#$, that is, a collection of subsets $A^\#(\sigma) \subseteq 
A(\sigma)$\/ closed under the application and containing the combinators ${\bf k}, {\bf s}$. There is a corresponding 
relative version of the extensional collapse. Longley's third theorem then states that when $(\cA,\cA^\#)$\/ is a
relative TPCA with $\cA$\/ full continuous and $\cA^\#$\/ effective, $\EC{\cA,\cA^\#}=\mathsf{RC}$, the total recursive
continuous functionals. The latter may be taken to be $\EC{\mathcal{K}_2,\mathcal{K}_2^{\rm eff}}$.

\subsection{Higher-order computation in toposes}

Since toposes are cartesian closed we can also consider higher type computability in toposes. Let us consider this first
in the case of the effective topos. The following result already appears in ~\cite{HylandJ:efft}:

\begin{theorem}[Hyland~\cite{HylandJ:efft}]
The total functionals of higher type over the NNO in $\eff$\/ are precisely the hereditarily effective operations.
\end{theorem}

Next, consider the \emph{Mulry topos}; this is the topos of sheaves on the monoid of total computable functions,
with the canonical topology. (The latter amounts to taking as basic coverings sets $\{f_1, \ldots, f_k\}$\/ for
which $\bigcup_{i=1}^kIm(f_i)=\NN$.) For the following result, a functional $G:\NN^\NN \to \NN$\/ is 
called \emph{Banach-Mazur} when for each computable $h:\NN^2 \to \NN$, the composite $G \circ \tilde{h}$\/ is
computable, where $\tilde{h}:\NN \to \NN^\NN$\/ is the transpose of $h$. 

\begin{theorem}[Mulry~\cite{MulryP:genbmf}]
The functionals $\csynt{N}^\csynt{N} \to \csynt{N}$\/ in the Mulry topos are precisely the Banach-Mazur functionals. 
\end{theorem}

Finally, let us consider the free topos. What are the total functionals of pure type in the free topos, i.e. arrows 
$\csynt{N}^{(k)}\rightarrow \csynt{N}^{(\ell)}$, $k,\ell \geq 1$?  
This question is answered in an interesting paper of A. Scedrov~\cite{Sced88}.

\begin{theorem}[Scedrov~\cite{Sced88}]  Let $\cF$ be the free topos, let $\cC$ be the free CCC with NNO
and let $\cF_{\cC}$ be the full subcategory of $\cF$ generated by $\cC$. 
The morphisms of $\cF_{\cC}$ are precisely those Kleene computable functionals that are 
provably total in the internal logic of $\cF$.  
\end{theorem}

The proof uses a gluing (or Friedman Realizability) argument (cf.~\cite{LS80,LS86,Wraith74}) together
with induction.  The ambient set theory is the free topos $\cF$ itself, and, for each type level $j\geq 0$ we construct 
an Effective topos $\eff(j)$ internally in the free topos, gluing it along a certain left exact functor 
$\Delta: \eff(j)\rightarrow \cF$.

 \subsection{Complexity Theory}
While classical computability theory is often concerned with the degree of unsolvability of various problems,
the branch most relevant to computer science is that of complexity theory, where one classifies solvable problems
according to the time and/or resources their solutions require. In particular, one is interested in
\emph{complexity classes} and the connections between those. For example, the class PTIME consists of
problems whose solution (regarded as a function of the input value $n \in \NN$) requires $p(n)$\/ steps
(of a deterministic Turing machine, say), where $p$\/ is a polynomial with positive integer coefficients. We refer
to~\cite{HarelD:alg} for an introduction.
 
\medskip 
 Early work in Implicit Computational Complexity by Martin Hofmann 
e.g.~\cite{Hof00} used complexity-bounded combinatory algebras 
and realizability to study logics of bounded complexity.  A BCK algebra is an applicative structure
$\cA$\/ having the combinators $\bbb, \ccc, \kkk$, where (still associating to the left)
\[ \bbb xyz=x(yz) \; ; \qquad \ccc xyz=xzy \;  ; \qquad \kkk xy=x.\]
Any total PCA is a BCK algebra, but not vice versa: the diagonal $x \mapsto xx$\/ is generally not computable in 
a BCK algebra. One of the results
in ~\cite{Hof00} shows that there is a BCK algebra structure on the natural numbers capturing PTIME computation:

\begin{theorem}[Hofmann~\cite{Hof00}]
There exists a BCK algebra structure on $\NN$\/ such that the computable maps w.r.t. this structure are precisely the
polynomial-time computable functions.
\end{theorem}

Related applications of such bounded combinatory algebras 
(to reprove the theorem that the 
  representable   functions of Bounded Linear Logic are exactly those in PTIME) appear in~\cite{HofScott04}.

   Recent work in Turing categories has focussed on the following general question:
   which complexity classes (e.g. LINEAR, PTIME, LOGSPACE, etc.) can occur
   as the total maps of a Turing Category?  Of course, such a Turing category cannot be
   a subcategory of $\Par$, since it would then necessarily contain all total computable functions. 
   Hence, it has to be a category whose global sections functor is not faithful. 
   
  The paper~\cite{CockettJ:T6} explores the area in more detail.
  Their main theorem characterizes when a Cartesian Category $\cC$ with a Universal
   Object $U$, a pair of disjoint elements \{{\sf t}, {\sf f}\}, and various abstract coding 
   retract structure can arise as the total maps of a Turing Category. The construction
   makes use of the idea that the given retract structure allows one to simulate a simple
   programming language. Passing to the presheaf topos of $\cC$ then provides the 
   required structure of a trace (on the coproduct) for implementing this language to obtain a PCA.
    
   As a consequence of this characterization, one obtains the following corollary:
       
  \begin{corollary}  Any countable Cartesian category with 
  a universal object $U$ and a pair of disjoint elements is the total maps of a Turing category.
  \end{corollary}

In order to apply this result to show that a particular complexity class arises as the total maps of a Turing category,
one is thus required to establish that the class in question admits the required closure conditions and pairing operations. 
For example, the classes  of LINEAR and PTIME maps (between binary numbers) can 
be shown to meet these requirements~\cite{CockettJ:T6}.
However, it is not fully understood for which complexity classes this is possible. 

\section*{Conclusion}  
We hope that we have shown in this -admittedly biased- overview of categorical recursion theory 
how various of Lambek's seminal ideas have initiated and inspired numerous strands of research that are still being pursued today.

We also hope to have conveyed to the reader that there are still many interesting unanswered questions and 
relatively unexplored facets of categorical recursion theory that deserve further investigation.

\bibliographystyle{plain}

\end{document}